 \documentclass[final,3p,times]{elsarticle}

\usepackage{amssymb}
 \usepackage{amsthm}
 \usepackage[T1]{fontenc}
 
\usepackage[ansinew]{inputenc} 
\usepackage{mathrsfs}
\usepackage{euscript}
\usepackage{natbib}

\newtheorem{theorem}{Theorem}
\newtheorem{corollary}{Corollary}
\newtheorem{definition}{Definition}
\newtheorem{lemma}{Lemma}
\newtheorem{proposition}{Proposition}
\newdefinition{rmk}{ Remark}
\newproof{pf}{Proof}

\begin{document}

\begin{frontmatter}



\title{ The Cauchy problem for heat equation with fractional Laplacian and exponential nonlinearity}

 \author[A.Fino]{A. Z. FINO}
 \address[A.Fino]{LaMA-Liban, Lebanese University, Faculty of Sciences, Department of Mathematics, P.O. Box 826 Tripoli, Lebanon}


 \ead{ahmad.fino01@gmail.com; afino@ul.edu.lb}

 \author[K1,K2,K3]{M. KIRANE}
\address[K1]{LaSIE, P\^ole Sciences et Technologies, Universit\'e de La Rochelle, Avenue Michel Cr\'epeau, 17031 La Rochelle, France}
\address[K2]{NAAM Research Group, Department of Mathematics,
King Abdulaziz University, P.O. Box 80203,
Jeddah 21589, Saudi Arabia}
\address[K3]{RUDN University, 6 Miklukho-Maklay St, Moscow 117198, Russia}
 \ead{mokhtar.kirane@univ-lr.fr}

 \begin{abstract}
We consider the Cauchy problem for heat equation with fractional Laplacian and exponential nonlinearity. We establish local well-posedness result in Orlicz spaces. We derive the existence of global solutions for small initial data. We obtain decay estimates for large time in Lebesgue spaces. 
\end{abstract}


\end{frontmatter}

\section{Introduction}

This paper concerns the Cauchy problem for the following  heat equation 
\begin{equation}\label{eq1}
\left\{\begin{array}{ll}
\,\, \displaystyle {u_{t}+(-\Delta)^{\beta/2}u =f(u),} &\displaystyle {t>0,x\in {\mathbb{R}^n},}\\
{}\\
\displaystyle{u(0,x)=  u_0(x),\qquad\qquad}&\displaystyle{x\in {\mathbb{R}^n},}
\end{array}
\right.
\end{equation} 
where $u$ is a real-valued unknown function, $0<\beta\leq 2$, $n\geq 1$, and $f:\mathbb{R}\rightarrow \mathbb{R}$ having an exponential growth at infinity ($f(u)\sim e^{|u|^p}$, $p>1$, for large $u$) with $f(0)=0$. Hereafter, $\|\cdotp\|_q$ $(1\leq q\leq\infty)$ stands for the usual $L^q(\mathbb{R}^n)$-norm.\\

When $f(u)=|u|^{p-1}u$, the Lebesgue spaces are adapted to study our problem (cf. \cite{Brezis,weissler1,weissler2,weissler3}). By analogy, we consider  the Orlicz spaces \cite{Orlicz} in order to study heat equations with exponential nonlinearities. 
The Orlicz space 
$$\textnormal{exp}\, L^p(\mathbb{R}^n)=\left\{u\in L^1_{loc}(\mathbb{R}^n);\,\int_{\mathbb{R}^n}\left(\textnormal{exp}\left(\frac{|u(x)|^p}{\lambda^p}\right)-1\right)\,dx<\infty,\, \textnormal{for some}\, \lambda>0\right\},$$
endowed with the Luxemburg norm
$$\|u\|_{\textnormal{exp~$L^p(\mathbb{R}^n)$}}:=\inf\left\{\lambda>0;\, \int_{\mathbb{R}^n}\left(\textnormal{exp}\left(\frac{|u(x)|^p}{\lambda^p}\right)-1\right)\,dx\leq 1\right\}$$
 is a Banach space.
For the local well-posedness we use the space
\begin{eqnarray*}\textnormal{exp}\, L^p_0(\mathbb{R}^n)=\Big\{u\in \textnormal{exp}\, L^p(\mathbb{R}^n); \textnormal{there exists $\{u_n\}_{n=1}^\infty\subset C^\infty_0(\mathbb{R}^n)$}\qquad\\
\quad\textnormal{such that $\lim_{n\rightarrow\infty}\|u_n-u\|_{\textnormal{exp~$L^p(\mathbb{R}^n)$}}=0$}\Big\}.
\end{eqnarray*}
It is also know (see Ioku, Ruf, and Terraneo \cite{Ioku}, Majdoub et al. \cite{Majdoub1,Majdoub2}) that
$$\textnormal{exp}\, L^p_0(\mathbb{R}^n)=\left\{u\in L^1_{loc}(\mathbb{R}^n);\,\int_{\mathbb{R}^n}\left(\textnormal{exp}\left(\alpha|u(x)|^p\right)-1\right)\,dx<\infty,\, \textnormal{for every}\, \alpha>0\right\}.$$

When $\beta=2$ (i.e. the standard heat equation) and $p=2$, Ioku \cite{Ioku2011} proved the existence of global solutions in $\textnormal{exp}\, L^2(\mathbb{R}^n)$ of (\ref{eq1}) under the condition (\ref{eq4}) below with $m=1+\frac{4}{n}$. Later,  Ioku et al. \cite{Ioku} studied the local nonexistence  of solutions of (\ref{eq1}) for certain data in $\textnormal{exp}\, L^2(\mathbb{R}^2)$, and the well-posedness of (\ref{eq1}) in the subspace $\textnormal{exp}\, L^2_0(\mathbb{R}^2)$ under the condition (\ref{eq3}) below. In \cite{Furioli}, Furioli et al.  considered the asymptotic behavior and decay estimates of the global solutions of (\ref{eq1}) in $\textnormal{exp}\, L^2(\mathbb{R}^n)$ when $f(u)=|u|^{4/n}ue^{u^2}$. Next, Majdoub et al. \cite{Majdoub1} proved the local well-posedness in $\textnormal{exp}\, L^2_0(\mathbb{R}^n)$ (if $f$ satisfies (\ref{eq3}) below with $m\geq 1+\frac{8}{n}$) and the global existence under small initial data in $\textnormal{exp}\, L^2(\mathbb{R}^n)$ (if $f$ satisfies (\ref{eq4}) below) for the biharmonic heat equation (i.e. $u_t+\Delta^2 u=f(u)$). Finally,  when $\beta=2$, $p>1$ and $m\geq 1+\frac{2p}{n}$, Majdoub and Tayachi \cite{Majdoub2} proved not only the local well-posedness in $\textnormal{exp}\, L^p_0(\mathbb{R}^n)$ but also the global existence of solutions, when $\frac{n(p-1)}{2}>p$, under small initial data in $\textnormal{exp}\, L^p(\mathbb{R}^n)$ of (\ref{eq1}) and analyzed their decay estimates, while the case of $\frac{n(p-1)}{2}\leq p$ is recently completed in \cite{Majdoub3}. In this paper, we generalize the papers of  \cite{Majdoub2, Majdoub3} to the fractional laplacian case.  \\

In order to state our main results, we note that the linear semigroup $e^{-t(-\Delta)^{\beta/2}}$ is continuous at $t=0$ in exp~$L^p_0(\mathbb{R}^n)$ (see Proposition \ref{prop2}) which is not the case in exp~$L^p(\mathbb{R}^n)$ (cf. \cite{Ioku} in the case of $\beta=2$), therefore, we have to define two kinds of mild solutions, the standard one where the space exp~$L^p_0(\mathbb{R}^n)$ is used, and the weak-mild solution where we use the space exp~$L^p(\mathbb{R}^n)$.
\begin{definition}(Mild solution)\label{mild}\\
Given $u_0\in \textnormal{exp}\, L^p_0(\mathbb{R}^n)$ and $T>0$. We say that $u$ is a mild solution for the Cauchy problem (\ref{eq1}) if $u\in C([0,T];\textnormal{exp}\, L^p_0(\mathbb{R}^n))$ satisfying
\begin{equation}\label{eq2}
u(t)=e^{-t(-\Delta)^{\beta/2}}u_0+\int_0^te^{-(t-s)(-\Delta)^{\beta/2}}f(u(s))\,ds,
\end{equation}
where $e^{-t(-\Delta)^{\beta/2}}$ is defined in (\ref{homogenous}) below.
\end{definition}
\begin{definition}(Weak-mild solution)\\
Given $u_0\in \textnormal{exp}\, L^p(\mathbb{R}^n)$ and $T>0$.  We say that $u$ is a weak-mild solution for the Cauchy problem (\ref{eq1}) if $u\in L^\infty((0,T);\textnormal{exp}\, L^p(\mathbb{R}^n))$ satisfying the associated integral equation (\ref{eq2}) in $\textnormal{exp}\, L^p(\mathbb{R}^n)$ for almost all $t\in(0,T)$ and $u(t)\rightarrow u_0$ in the weak$^{*}$ topology as $t\rightarrow 0$.
\end{definition}
We recall that $u(t)\rightarrow u_0$ in weak$^{*}$ sense if and only if 
$$\lim_{t\rightarrow 0}\int_{\mathbb{R}^n}[u(t,x)\varphi(x)-u_0(x)\varphi(x)]\,dx=0,\quad\textnormal{for every}\, \varphi\in L^1(\ln L)^{1/p}(\mathbb{R}^n),$$
where
$$L^1(\ln L)^{1/p}(\mathbb{R}^n):=\left\{f\in L^1_{loc}(\mathbb{R}^n);\,\int_{\mathbb{R}^n}|f(x)|\ln^{1/p}(2+|f(x)|)\,dx<\infty\right\}$$
 is a predual space of $\textnormal{exp}\, L^p(\mathbb{R}^n)$ (see \cite{Bennett,Rao}).\\ 
 
 First, we interest in the local well-posedness. We assume that $f$ satisfies
 \begin{equation}\label{eq3}
f(0)=0,\qquad |f(u)-f(v)|\leq C|u-v|(e^{\lambda |u|^p}+e^{\lambda |v|^p}),\quad \hbox{for all}\,\, u,v\in\mathbb{R},
\end{equation}
for some constants $C>0$, $p>1$, and $\lambda>0$. Typical example satisfying (\ref{eq3}) is: $f(u)=\pm ue^{|u|^p}$.
\begin{theorem}\label{theo1}(Local well-posedness)\\
Let $n\geq1$, $p>1$, and $0<\beta< 2$. Suppose that $f$ satisfies $(\ref{eq3})$. Given $u_0\in \textnormal{exp}\, L^p_0(\mathbb{R}^n)$, there exist a time $T=T(u_0)>0$ and a unique mild solution $u\in C([0,T];\textnormal{exp}\, L^p_0(\mathbb{R}^n))$ to $(\ref{eq1})$.
\end{theorem}

Next, our second interest is the global existence and the decay estimate. In this case, the behaviour of $f(u)$ near $u=0$ plays a crucial role, therefore the following behaviour near zero will be allowed
$$|f(u)|\sim|u|^m,$$
where $\frac{n(m-1)}{\beta}\geq p$. More precisely, we suppose that
 \begin{equation}\label{eq4}
f(0)=0,\qquad |f(u)-f(v)|\leq C|u-v|(|u|^{m-1}e^{\lambda |u|^p}+|v|^{m-1}e^{\lambda |v|^p}),\quad \hbox{for all}\,\, u,v\in\mathbb{R},
\end{equation}
where $\frac{n(m-1)}{\beta}\geq p>1$, $C>0$, and $\lambda>0$ are constants. Typical example satisfying (\ref{eq4}) is: $f(u)=\pm |u|^{m-1}ue^{|u|^p}$ where $m\geq 1+\frac{\beta p}{n}$; we note that the global existence in the case $m= 1+\frac{\beta p}{n}$ is presented in \cite[Section~8]{Furioli} without any proof.
\begin{theorem}\label{theo2}(Global existence)\\
Let $n\geq1$, $1<p\leq\frac{n(m-1)}{\beta}$, and $0<\beta< 2$. Suppose that $f$ satisfies $(\ref{eq4})$ for $m\geq p$. Then there exists a positive constant $\varepsilon>0$ such that every initial data $u_0\in \textnormal{exp}\, L^p(\mathbb{R}^n)$ with $\|u_0\|_{\textnormal{exp~$L^p(\mathbb{R}^n)$}}\leq\varepsilon$, there exists  a global weak-mild solution $u\in L^\infty((0,\infty);\textnormal{exp}\, L^p(\mathbb{R}^n))$ to $(\ref{eq1})$ satisfying
 \begin{equation}\label{eq5}
\lim_{t\rightarrow0}\left\|u(t)-e^{-t(-\Delta)^{\beta/2}}u_0\right\|_{\textnormal{exp~$L^p(\mathbb{R}^n)$}}=0.
\end{equation}
Moreover, there exists a constant $C>0$ such that
 \begin{equation}\label{eq35}
 \|u(t)\|_{L^q(\mathbb{R}^n)} \leq  C t^{-\sigma}, \quad\hbox{for all}\,\,\,t>0,
 \end{equation}
where
$$ \sigma=\frac{1}{m-1}-\frac{n}{\beta q}>0,$$
and
$$
\left\{\begin{array}{ll}
\,\, \displaystyle {\frac{n(m-1)}{\beta}<q< \frac{n(m-1)}{\beta}\frac{1}{(2-m)_+}}, &\displaystyle {\hbox{if}\,\,\,\beta=\frac{n(p-1)}{p},}\\
{}\\
\displaystyle{\frac{n(m-1)}{\beta}<q<\frac{n(m-1)}{\beta}\frac{1}{(2-m)_+},\qquad}&\displaystyle{\hbox{if}\,\,\,\beta<\frac{n(p-1)}{p},}\\
{}\\
\displaystyle{\frac{(m-1)p}{p-1}<q<\frac{n(m-1)}{\beta}\frac{1}{(2-m)_+},\qquad}&\displaystyle{\hbox{if}\,\,\,\beta>\frac{n(p-1)}{p} \,\&\,(2-m)_+<\frac{n(p-1)}{p\beta},}\\
\end{array}
\right.
$$
with $(\cdotp)_+$ stands for the positive part.
\end{theorem}
\begin{rmk}
In Theorem \ref{theo2}, we have to distinguish 3 cases: $\beta<\frac{n(p-1)}{p}$, $\beta>\frac{n(p-1)}{p}$, and $\beta=\frac{n(p-1)}{p}$. We note that in the case of $\beta>\frac{n(p-1)}{p}$ we have to take $m>p$. Indeed, if $m=p$, it follows that $\beta>\frac{n(m-1)}{m}$, but $n(m-1)/\beta\geq p$, which implies that $\beta\leq\frac{n(m-1)}{m}$, therefore $\frac{n(p-1)}{p}<\frac{n(p-1)}{p}$; contradiction.
\end{rmk}

This paper is organized as follows: in Section \ref{sec2}, we present several preliminaries. Section \ref{sec3} contains the proof of the local well-posedness theorem (Theorem \ref{theo1}). Finally, we prove the global existence theorem (Theorem \ref{theo2}) in Section \ref{sec4}.


\section{Preliminaries}\label{sec2}
\subsection{Orlicz spaces: basic properties}
In this section we present the definition of the so-called Orlicz spaces on $\mathbb{R}^n$ and some related properties.  More details and complete presentations can be found in \cite{Adams,Rao,Trudinger}.
\begin{definition}(Orlicz space)\\
Let $\phi:\mathbb{R}^+\rightarrow \mathbb{R}^+$ be a convex increasing function such that
$$\phi(0)=0=\lim_{s\rightarrow 0^+}\phi(s), \quad \lim_{s\rightarrow \infty}\phi(s)=\infty.$$
The Orlicz space $L^\phi(\mathbb{R}^n)$ is defined by
$$L^\phi(\mathbb{R}^n)=\left\{u\in L^1_{loc}(\mathbb{R}^n);\,\int_{\mathbb{R}^n}\phi\left(\frac{|u(x)|}{\lambda}\right)\,dx<\infty,\,\, \textnormal{for some}\, \lambda>0\right\},$$
endowed with the Luxemburg norm
$$\|u\|_{L^\phi(\mathbb{R}^n)}:=\inf\left\{\lambda>0;\,\int_{\mathbb{R}^n}\phi\left(\frac{|u(x)|}{\lambda}\right)\,dx\leq 1\right\}.$$
\end{definition}
On the other hand, we denote by
$$L^\phi_0(\mathbb{R}^n)=\left\{u\in L^1_{loc}(\mathbb{R}^n);\,\int_{\mathbb{R}^n}\phi\left(\frac{|u(x)|}{\lambda}\right)\,dx<\infty,\,\, \textnormal{for every}\, \lambda>0\right\}.$$
It can be shown (as in Ioku et al. \cite{Ioku}) that
$$L^\phi_0(\mathbb{R}^n)=\overline{C^\infty_0(\mathbb{R}^n)}^{\|\cdotp\|_{L^\phi}}=\,\textnormal{the closure of $C^\infty_0(\mathbb{R}^n)$ in $L^\phi(\mathbb{R}^n)$}.$$
It is known that $(L^\phi(\mathbb{R}^n),\|\cdotp\|_{L^\phi(\mathbb{R}^n)})$ and $(L^\phi_0(\mathbb{R}^n),\|\cdotp\|_{L^\phi(\mathbb{R}^n)})$ are Banach spaces. Note that, if $\phi(s)=s^p$, $1\leq p<\infty$, then $L^\phi(\mathbb{R}^n)=L^\phi_0(\mathbb{R}^n)=L^p(\mathbb{R}^n)$, and if $\phi(s)=e^{s^p}-1$, $1\leq p<\infty$, then $L^\phi(\mathbb{R}^n)$ is the space $\textnormal{exp}\, L^p(\mathbb{R}^n)$, while $L^\phi_0(\mathbb{R}^n)$ is the space $\textnormal{exp}\, L^p_0(\mathbb{R}^n)$. Moreover, for $u\in L^\phi$ and $K:=\|u\|_{L^\phi(\mathbb{R}^n)}>0$, we can easy check by the definition of the infimum that
$$\left\{\lambda>0,\,\,\int_{\mathbb{R}^n}\phi\left(\frac{|u(x)|}{\lambda}\right)\,dx\leq 1\right\}=[K;\infty[,$$
in particular
\begin{equation}\label{eq9}
\int_{\mathbb{R}^n}\phi\left(\frac{|u(x)|}{\|u\|_{L^\phi(\mathbb{R}^n)}}\right)\,dx\leq 1.
\end{equation}
The following Lemmas summarize the embedding between Orlicz  and Lebesgue spaces.
\begin{lemma}\cite[Lemma 2.3]{Majdoub2}\label{lemma1}\\
For every $1\leq q\leq p$, we have $L^q(\mathbb{R}^n)\cap L^\infty(\mathbb{R}^n)\hookrightarrow \textnormal{exp}\, L^p_0(\mathbb{R}^n)\hookrightarrow \textnormal{exp}\, L^p(\mathbb{R}^n)$, more precisely
\begin{equation}\label{eq6}
\|u\|_{\textnormal{exp~$L^p(\mathbb{R}^n)$}}\leq\frac{1}{(\ln 2)^{1/p}}(\|u\|_q+\|u\|_\infty).
\end{equation}
\end{lemma}
\noindent Similarly, we have
\begin{lemma}\label{lemma2}${}$\\
Let $\phi(s)=e^{s^p}-1-s^p$, $p>1$. For every $ q\leq 2p$, we have $L^q(\mathbb{R}^n)\cap L^\infty(\mathbb{R}^n)\hookrightarrow L^\phi_0(\mathbb{R}^n)\hookrightarrow L^\phi(\mathbb{R}^n)$, more precisely
\begin{equation}\label{eq7}
\|u\|_{ L^\phi(\mathbb{R}^n)}\leq C(p)(\|u\|_q+\|u\|_\infty).
\end{equation}
\end{lemma}
\proof Let $g(s)=e^{s^p}-s^p$; $g$ is a strictly increasing. Let $\alpha\geq C(p)(\|u\|_q+\|u\|_\infty)$ where $C(p):=1/g^{-1}(2)$, then 
\begin{eqnarray*}
&{}&\int_{\mathbb{R}^n}\left(\textnormal{exp}\left(\frac{|u(x)|^p}{\alpha^p}\right)-1-\left(\frac{|u(x)|^p}{\alpha^p}\right)\right)\,dx\\
&{}&\,=\sum_{k=2}^{\infty}\frac{1}{k!\alpha^{pk}}\|u\|_{pk}^{pk}\\
&{}&\,\leq \sum_{k=2}^{\infty}\frac{1}{k!\alpha^{pk}}(\|u\|_{q}+\|u\|_\infty)^{pk}\\
&{}&\,=\textnormal{exp}\left(\frac{\|u\|_{q}+\|u\|_\infty}{\alpha}\right)^p-1-\left(\frac{\|u\|_{q}+\|u\|_\infty}{\alpha}\right)^p\\
&{}&\,=g\left(\frac{\|u\|_{q}+\|u\|_\infty}{\alpha}\right)-1\\
&{}&\,\leq 1,
\end{eqnarray*}
where we have used the interpolation inequality $\|u\|_r\leq \|u\|_{q}^{q/r}\|u\|_{\infty}^{1-q/r}\leq \|u\|_{q}+\|u\|_\infty$ for all $q\leq r\leq \infty$ and all $u\in L^q\cap L^\infty$. Therefore
$$[C(p)(\|u\|_q+\|u\|_\infty);\infty[\subseteq \left\{\alpha>0;\,\int_{\mathbb{R}^n}\phi\left(\frac{|u(x)|}{\alpha}\right)\,dx\leq 1\right\},$$
which implies that 
\begin{eqnarray*}
\|u\|_{ L^\phi(\mathbb{R}^n)}&=&\inf \left\{\alpha>0;\,\int_{\mathbb{R}^n}\phi\left(\frac{|u(x)|}{\alpha}\right)\,dx
\leq 1\right\}\\
&\leq& \inf \left\{\alpha>0;\,\alpha\in [C(p)(\|u\|_q+\|u\|_\infty);\infty[\right\}\\
&=& C(p)(\|u\|_q+\|u\|_\infty).
\end{eqnarray*}
\hfill$\square$
\begin{lemma}\cite[Lemma 2.4]{Majdoub2}\label{lemma3}\\
For every $1\leq p\leq q<\infty$, we have $\textnormal{exp}\, L^p(\mathbb{R}^n)\hookrightarrow L^q(\mathbb{R}^n)$, more precisely
\begin{equation}\label{eq8}
\|u\|_q\leq\left(\Gamma\left(\frac{q}{p}+1\right)\right)^{1/q} \|u\|_{\textnormal{exp~$L^p(\mathbb{R}^n)$}},
\end{equation}
where $\Gamma$ is the gamma function.
\end{lemma}
Next, we present some definitions and results concerning the
fractional Laplacian that will be used hereafter. The fundamental solution $S_\beta$ of the usual linear fractional diffusion equation
\begin{equation}\label{homo}
u_t+(-\Delta)^{\beta/2}u=0,\quad \beta\in(0,2],\;
x\in\mathbb{R}^n,\;t>0,
\end{equation}
can be represented via the Fourier transform by
\begin{equation}\label{FS}
S_\beta(t)(x):=S_\beta(x,t)=\frac{1}{(2\pi)^{n/2}}\int_{\mathbb{R}^n}e^{ix.\xi-t|{\xi}|^\beta}\,d\xi.
\end{equation}
This mean that the solution of (\ref{homo}) with any initial data $u(0)=u_0$ can be written as
\begin{equation}\label{homogenous}
u(x,t)=S_\beta(x,t)\ast u_0 (x)=:e^{-t(-\Delta)^{\beta/2}}u_0,
\end{equation}
where $e^{-t(-\Delta)^{\beta/2}}$ is a  strongly continuous semigroup on $L^p(\mathbb{R}^n)$, $p>1$, generated by the fractional power $-(-\Delta)^{\beta/2}$. Moreover, $S_\beta$ satisfies
\begin{equation}\label{P_1}
    S_\beta(1)\in L^\infty(\mathbb{R}^n)\cap
L^1(\mathbb{R}^n),\quad
S_\beta(x,t)\geq0,\quad\int_{\mathbb{R}^n}S_\beta(x,t)\,dx=1,
\end{equation}
\noindent for all $x\in\mathbb{R}^n$ and $t>0.$ Hence, using Young's inequality for the convolution
and the following self-similar
form $S_\beta(x,t)=t^{-n/\beta}S_\beta(xt^{-1/\beta},1),$ we get the $L^r-L^q$ estimate
\begin{equation}\label{P_2}
\left\|e^{-t(-\Delta)^{\beta/2}} v\right\|_q\;\leq \;C\,t^{-\frac{n}{\beta}(\frac{1}{r}-\frac{1}{q})}\|v\|_r,
\end{equation}
for all $v\in L^r(\mathbb{R}^n)$ and all $1\leq r\leq q\leq\infty,$ $t>0$, where $C>0$ is a positive constant depending only on $n$. In particular, using Young's inequality for the convolution and (\ref{P_1}), we have
\begin{equation}\label{P_3}
\left\|e^{-t(-\Delta)^{\beta/2}} v\right\|_q=\left\|S_\beta(x,t)\ast v\right\|_q\;\leq \left\|S_\beta(t)\right\|_1\|v\|_q=\|v\|_q,
\end{equation}
for all $v\in L^q(\mathbb{R}^n)$ and all $1\leq q\leq\infty,$ $t>0.$ \\
The following proposition is a generalization of Proposition 3.2 in \cite{Majdoub2} and it is presented (without proof) by Furioli et al. \cite[Lemma 3.1]{Furioli}.
\begin{proposition}\label{prop1}${}$\\
Let $1\leq q\leq p$, $1\leq r\leq\infty$, and $0<\beta\leq 2$. Then the following estimates hold.
\begin{enumerate}[(i)]
\item $\left\|e^{-t(-\Delta)^{\beta/2}} \varphi\right\|_{\textnormal{exp~$L^p(\mathbb{R}^n)$}}\leq \left\|\varphi\right\|_{\textnormal{exp~$L^p(\mathbb{R}^n)$}}$, \,\,\, for all $t>0,$\, $\varphi\in {\textnormal{exp~$L^p(\mathbb{R}^n)$}}$.
\item $\left\|e^{-t(-\Delta)^{\beta/2}} \varphi\right\|_{\textnormal{exp~$L^p(\mathbb{R}^n)$}}\leq C\,t^{-\frac{n}{\beta q}}\left(\ln( t^{-\frac{n}{\beta}}+1)\right)^{-1/p}\left\|\varphi\right\|_{q}$, \,\,\, for all $t>0,$\, $\varphi\in L^q(\mathbb{R}^n)$.
\item $\left\|e^{-t(-\Delta)^{\beta/2}} \varphi\right\|_{\textnormal{exp~$L^p(\mathbb{R}^n)$}}\leq \frac{1}{(\ln 2)^{1/p}}\left[C\,t^{-\frac{n}{\beta r}}\left\|\varphi\right\|_{r}+\left\|\varphi\right\|_{q}\right]$, \,\,\, for all $t>0,$\, $\varphi\in L^r(\mathbb{R}^n)\cap L^q(\mathbb{R}^n)$.
\end{enumerate}
\end{proposition}
\proof We start by proving (i). For any $\lambda>0$, using (\ref{P_3}) and Taylor expansion, we have
\begin{eqnarray*}
\int_{\mathbb{R}^n}\left(\textnormal{exp}\left(\frac{|e^{-t(-\Delta)^{\beta/2}} \varphi|^p}{\lambda^p}\right)-1\right)\,dx&=&\sum_{k=1}^{\infty}\frac{\left\|e^{-t(-\Delta)^{\beta/2}} \varphi\right\|_{pk}^{pk}}{k!\lambda^{pk}}\\
&\leq& \sum_{k=1}^{\infty}\frac{\left\|\varphi\right\|_{pk}^{pk}}{k!\lambda^{pk}}\\
&=&\int_{\mathbb{R}^n}\left(\textnormal{exp}\left(\frac{|\varphi|^p}{\lambda^p}\right)-1\right)\,dx.
\end{eqnarray*}
Then
\begin{eqnarray*}
&{}&\left\{\lambda>0;\,\,\int_{\mathbb{R}^n}\left(\textnormal{exp}\left(\frac{|\varphi|^p}{\lambda^p}\right)-1\right)\,dx\leq 1\right\}\\
&{}&\subseteq \left\{\lambda>0;\,\,\int_{\mathbb{R}^n}\left(\textnormal{exp}\left(\frac{|e^{-t(-\Delta)^{\beta/2}} \varphi|^p}{\lambda^p}\right)-1\right)\,dx\leq 1\right\},
\end{eqnarray*}
and therefore
\begin{eqnarray*}
\left\|e^{-t(-\Delta)^{\beta/2}} \varphi\right\|_{\textnormal{exp~$L^p(\mathbb{R}^n)$}}&=&\inf\left\{\lambda>0;\,\,\int_{\mathbb{R}^n}\left(\textnormal{exp}\left(\frac{|e^{-t(-\Delta)^{\beta/2}} \varphi|^p}{\lambda^p}\right)-1\right)\,dx\leq 1\right\}\\
&\leq& \inf\left\{\lambda>0;\,\,\int_{\mathbb{R}^n}\left(\textnormal{exp}\left(\frac{|\varphi|^p}{\lambda^p}\right)-1\right)\,dx\leq 1\right\}\\
&=&\left\|\varphi\right\|_{\textnormal{exp~$L^p(\mathbb{R}^n)$}}.
\end{eqnarray*}
This proves (i). Similarly, to prove (ii), we use again (\ref{P_2}) and Taylor expansion. For any $\lambda>0$, we have
\begin{eqnarray*}
\int_{\mathbb{R}^n}\left(\textnormal{exp}\left(\frac{|e^{-t(-\Delta)^{\beta/2}} \varphi|^p}{\lambda^p}\right)-1\right)\,dx&=&\sum_{k=1}^{\infty}\frac{\left\|e^{-t(-\Delta)^{\beta/2}} \varphi\right\|_{pk}^{pk}}{k!\lambda^{pk}}\\
&\leq& \sum_{k=1}^{\infty}\frac{C^{pk}t^{-\frac{n}{\beta}(\frac{1}{q}-\frac{1}{pk})pk}\left\|\varphi\right\|_{q}^{pk}}{k!\lambda^{pk}}\\
&=&t^{\frac{n}{\beta}}\left(\textnormal{exp}\left(\frac{C t^{-\frac{n}{\beta q}}\|\varphi\|_q}{\lambda}\right)^p-1\right).
\end{eqnarray*}
As
$$t^{\frac{n}{\beta}}\left(\textnormal{exp}\left(\frac{C t^{-\frac{n}{\beta q}}\|\varphi\|_q}{\lambda}\right)^p-1\right)\leq 1\Longleftrightarrow \lambda\geq C\,t^{-\frac{n}{\beta q}}\left(\ln( t^{-\frac{n}{\beta}}+1)\right)^{-1/p}\left\|\varphi\right\|_{q},$$
we conclude that 
\begin{eqnarray*}
&{}&\left\{\lambda>0,\,\, \lambda\in [ C\,t^{-\frac{n}{\beta q}}\left(\ln( t^{-\frac{n}{\beta}}+1)\right)^{-1/p}\left\|\varphi\right\|_{q};\infty[\right\}\\
&{}&\subseteq \left\{\lambda>0,\,\, \int_{\mathbb{R}^n}\left(\textnormal{exp}\left(\frac{|e^{-t(-\Delta)^{\beta/2}} \varphi|^p}{\lambda^p}\right)-1\right)\,dx\leq 1\right\};
\end{eqnarray*}
whereupon
\begin{eqnarray*}
\left\|e^{-t(-\Delta)^{\beta/2}} \varphi\right\|_{\textnormal{exp~$L^p(\mathbb{R}^n)$}}&=&\inf\left\{\lambda>0,\,\, \int_{\mathbb{R}^n}\left(\textnormal{exp}\left(\frac{|e^{-t(-\Delta)^{\beta/2}} \varphi|^p}{\lambda^p}\right)-1\right)\,dx\leq 1\right\}\\
&\leq&\inf \left\{\lambda>0,\,\, \lambda\in [ C\,t^{-\frac{n}{\beta q}}\left(\ln( t^{-\frac{n}{\beta}}+1)\right)^{-1/p}\left\|\varphi\right\|_{q};\infty[\right\}\\
&=&C\,t^{-\frac{n}{\beta q}}\left(\ln( t^{-\frac{n}{\beta}}+1)\right)^{-1/p}\left\|\varphi\right\|_{q}.
\end{eqnarray*}
This proves (ii). Finally, to prove (iii), we use the embedding $L^q(\mathbb{R}^n)\cap L^\infty(\mathbb{R}^n)\hookrightarrow \textnormal{exp}\, L^p_0(\mathbb{R}^n)$ (\ref{eq6}); we get
$$\left\|e^{-t(-\Delta)^{\beta/2}} \varphi\right\|_{\textnormal{exp~$L^p(\mathbb{R}^n)$}}\leq\frac{1}{(\ln 2)^{1/p}}\left(\left\|e^{-t(-\Delta)^{\beta/2}} \varphi\right\|_q+\left\|e^{-t(-\Delta)^{\beta/2}} \varphi\right\|_\infty\right).$$
Using the $L^r-L^\infty$ and $L^q-L^q$ estimates (\ref{P_2}), we conclude that
$$\left\|e^{-t(-\Delta)^{\beta/2}} \varphi\right\|_{\textnormal{exp~$L^p(\mathbb{R}^n)$}}\leq\frac{1}{(\ln 2)^{1/p}}\left(\left\| \varphi\right\|_q+C\,t^{-\frac{n}{\beta r}}\left\| \varphi\right\|_r\right).$$
\hfill$\square$\\
We will also need the following smoothing results.
\begin{proposition}\label{prop2}${}$\\
If $\varphi\in \textnormal{exp~$L^p_0(\mathbb{R}^n)$}$, then $e^{-t(-\Delta)^{\beta/2}}\varphi\in C([0,\infty);\textnormal{exp~$L^p_0(\mathbb{R}^n)$}).$
\end{proposition}
\proof The proof of this proposition follows the same one of \cite[Proposition~3.7]{Majdoub1} by making the appropriate modifications. To be self-contained, we will present it in details. Let $\varphi\in \textnormal{exp~$L^p_0(\mathbb{R}^n)$}$. By (i) of Proposition \ref{prop1} and the definition of $\textnormal{exp~$L^p_0(\mathbb{R}^n)$}$, we have $e^{-t(-\Delta)^{\beta/2}}\varphi\in \textnormal{exp~$L^p_0(\mathbb{R}^n)$}$ for every $t>0$. Thus, by the linearity of the semigroup  $e^{-t(-\Delta)^{\beta/2}}$, it remains to prove the continuity at $t=0$,
$$\lim_{t\rightarrow 0}\left\|e^{-t(-\Delta)^{\beta/2}} \varphi-\varphi\right\|_{\textnormal{exp~$L^p(\mathbb{R}^n)$}}=0.$$
Since $\varphi\in \textnormal{exp~$L^p_0(\mathbb{R}^n)$}$, there exists a sequence $(\varphi_n)_n\subseteq C^\infty_0(\mathbb{R}^n)$ such that\\
$\lim_{n\rightarrow \infty}\left\|\varphi_n-\varphi\right\|_{\textnormal{exp~$L^p$}}=0$. By (\ref{eq6}), and estimation (i) of Proposition \ref{prop1}, we obtain
\begin{eqnarray*}
&{}&\left\|e^{-t(-\Delta)^{\beta/2}} \varphi-\varphi\right\|_{\textnormal{exp~$L^p(\mathbb{R}^n)$}}\\
&{}&\,\leq\left\|e^{-t(-\Delta)^{\beta/2}} (\varphi-\varphi_n)\right\|_{\textnormal{exp~$L^p$}}+\left\|e^{-t(-\Delta)^{\beta/2}} \varphi_n-\varphi_n\right\|_{\textnormal{exp~$L^p$}}+\left\|\varphi_n-\varphi\right\|_{\textnormal{exp~$L^p$}}\\
&{}&\,\leq\frac{1}{(\ln 2)^{1/p}}\left(\|e^{-t(-\Delta)^{\beta/2}} \varphi_n-\varphi_n\|_p+\|e^{-t(-\Delta)^{\beta/2}} \varphi_n-\varphi_n\|_\infty\right)+2\left\|\varphi_n-\varphi\right\|_{\textnormal{exp~$L^p$}}.
\end{eqnarray*}
Since $\varphi_n\in C^\infty_0(\mathbb{R}^n)$, using the fact that $e^{-t(-\Delta)^{\beta/2}}$ is a strongly continuous semigroup on $L^r(\mathbb{R}^n)$ $(1<r\leq\infty)$, we have 
$$\lim_{t\rightarrow 0}\left(\|e^{-t(-\Delta)^{\beta/2}} \varphi_n-\varphi_n\|_p+\|e^{-t(-\Delta)^{\beta/2}} \varphi_n-\varphi_n\|_\infty\right)=0.$$ 
Hence
$$\limsup_{t\rightarrow 0}\left\|e^{-t(-\Delta)^{\beta/2}} \varphi-\varphi\right\|_{\textnormal{exp~$L^p(\mathbb{R}^n)$}}\leq 2\left\|\varphi_n-\varphi\right\|_{\textnormal{exp~$L^p$}},$$
for every $n\in\mathbb{N}$. This finishes the proof of the proposition. \hfill$\square$\\

It is known that $e^{-t(-\Delta)^{\beta/2}}$ is a $C^0$-semigroup on $L^p(\mathbb{R}^n)$. By Proposition \ref{prop2}, it is a $C^0$-semigroup on $\textnormal{exp~$L^p_0(\mathbb{R}^n)$}$.
\begin{lemma}\label{lemma10}\cite[Lemma 4.1.5]{CH}\\
Let $X$ be a Banach space and $g\in L^1(0,T;X)$, then $\displaystyle \int_0^te^{-(t-s)(-\Delta)^{\beta/2}}g(s)\,ds\in C([0,T];X)$. Moreover
$$\left\|\int_0^te^{-(t-s)(-\Delta)^{\beta/2}}g(s)\,ds\right\|_{L^\infty(0,T;X)}\leq \|g\|_{L^1(0,T;X)}.$$
\end{lemma}

The following lemmas are essential for the proof of the global existence (Theorem \ref{theo2}).
\begin{lemma}${}$\label{lemma6}\cite[Lemma~2.6]{Majdoub2}\\
Let $\lambda>0$, $1\leq q<\infty$ and $K>0$ be such that $\lambda q K^p\leq 1$. Assume that $u\in\textnormal{exp~$L^p(\mathbb{R}^n)$}$ satisfies
$$\|u\|_{\textnormal{exp~$L^p(\mathbb{R}^n)$}}\leq K,$$
then $\textnormal{exp}\left(\frac{|u|^p}{\lambda^p}\right)-1\in L^q(\mathbb{R}^n)$ and
$$\left\| e^{\lambda|u|^p}-1\right\|_{L^q(\mathbb{R}^n)}\leq \left(\lambda q K^p\right)^{1/q}.$$
\end{lemma}

\begin{lemma}${}$\label{lemma4}\\
Let $p>1$, $0<\beta\leq 2$ be such that $\beta<\frac{n(p-1)}{p}$. Then, for every $r>\frac{n}{\beta}$, there exists $C=C(n,p,\beta,r)$ such that
 $$\left\|\int_0^te^{-(t-s)(-\Delta)^{\beta/2}} g(s)\,ds\right\|_{L^\infty(0,\infty;\textnormal{exp~$L^p(\mathbb{R}^n)$})}\leq C\|g\|_{L^\infty(0,\infty;L^1(\mathbb{R}^n)\cap L^r(\mathbb{R}^n))},$$ 
for every $g\in L^\infty(0,\infty;L^1(\mathbb{R}^n)\cap L^r(\mathbb{R}^n))$.
\end{lemma}
\proof By Proposition \ref{prop1} (ii) with $q=1$, we have
\begin{equation}\label{eq10}
\left\|e^{-t(-\Delta)^{\beta/2}} \varphi\right\|_{\textnormal{exp~$L^p(\mathbb{R}^n)$}}\leq C\,t^{-\frac{n}{\beta }}\left(\ln( t^{-\frac{n}{\beta}}+1)\right)^{-1/p}\left\|\varphi\right\|_{1},
\end{equation}
for all $t>0$, $\varphi\in L^1(\mathbb{R}^n)\cap L^r(\mathbb{R}^n)$ ($\|\varphi\|_{L^1\cap L^r}=\|\varphi\|_{L^1}+\|\varphi\|_{L^r}$), while by Proposition \ref{prop1} (iii) with $q=1$, we obtain
\begin{equation}\label{eq11}
\left\|e^{-t(-\Delta)^{\beta/2}} \varphi\right\|_{\textnormal{exp~$L^p(\mathbb{R}^n)$}}\leq C\,(t^{-\frac{n}{\beta r}}+1)\left[\left\|\varphi\right\|_{r}+\left\|\varphi\right\|_{1}\right].
\end{equation}
Combining (\ref{eq10}) and (\ref{eq11}), we get
$$\left\|e^{-t(-\Delta)^{\beta/2}} \varphi\right\|_{\textnormal{exp~$L^p(\mathbb{R}^n)$}}\leq \kappa(t)\left[\left\|\varphi\right\|_{r}+\left\|g\right\|_{1}\right],$$
where 
$$\kappa(t)=\min\left\{C\,(t^{-\frac{n}{\beta r}}+1),C\,t^{-\frac{n}{\beta }}\left(\ln( t^{-\frac{n}{\beta}}+1)\right)^{-1/p}\right\}.$$
Due to the assumptions $\beta<\frac{n(p-1)}{p}$ and $r>\frac{n}{\beta}$, we see that $\kappa\in L^1(0,\infty)$. Thus, for $g\in L^\infty(0,\infty;L^1(\mathbb{R}^n)\cap L^r(\mathbb{R}^n))$, we have
\begin{eqnarray*}
\left\|\int_0^te^{-(t-s)(-\Delta)^{\beta/2}} g(s)\,ds\right\|_{\textnormal{exp~$L^p(\mathbb{R}^n)$}}&\leq& \int_0^t\left\|e^{-(t-s)(-\Delta)^{\beta/2}} g(s)\right\|_{\textnormal{exp~$L^p(\mathbb{R}^n)$}}\,ds\\
&\leq& \int_0^t\kappa(t-s)\left(\left\|g(s)\right\|_{L^1(\mathbb{R}^n)}+\left\|g(s)\right\|_{L^r(\mathbb{R}^n)}\right)\,ds\\
&\leq&\|g\|_{L^\infty(0,\infty;L^1(\mathbb{R}^n)\cap L^r(\mathbb{R}^n))} \int_0^\infty\kappa(s)\,ds,
\end{eqnarray*}
for every $t>0$. This proves Lemma \ref{lemma4}. \hfill$\square$\\

We remark that $\frac{n(p-1)}{p}$ may not included in $(0,2]$. So if $\frac{n(p-1)}{p}>2$, we have $\frac{n(p-1)}{p}>\beta$, and this case is recovered by Lemma \ref{lemma4}. If $\frac{n(p-1)}{p}\leq2$, we have three case to study: the case of $\beta<\frac{n(p-1)}{p}$ is done by Lemma \ref{lemma4}, and the case $\beta>\frac{n(p-1)}{p}$ can be done separately without using any kind of an a priori estimate, so it remains to study the case of $\beta=\frac{n(p-1)}{p}$ where we have a similar result as in Lemma \ref{lemma4}.  For this, we need to introduce an appropriate Orlicz space. Let $L^\phi(\mathbb{R}^n)$ this space, with $\phi(u)=e^{|u|^p}-1-|u|^p$, associated with its Luxemburg norm. From the definition of $\|\cdotp\|_{L^\phi}$, (\ref{eq8}), and the standard inequality $e^{\theta s}-1\leq\theta(e^s-1)$, $0\leq\theta\leq1$, $s\geq 0$, we can easily get
\begin{equation}\label{eq0}
C_1\|u\|_{\textnormal{exp~$L^p(\mathbb{R}^n)$}}\leq \|u\|_{L^p(\mathbb{R}^n)}+\|u\|_{L^\phi(\mathbb{R}^n)}  \leq C_2 \|u\|_{\textnormal{exp~$L^p(\mathbb{R}^n)$}},
\end{equation}
for some $C_1,C_2>0$. 
\begin{lemma}${}$\label{lemma5}\\
Let $p>1$, $0<\beta\leq 2$ be such that $\beta=\frac{n(p-1)}{p}$. Then, there exists $C=C(n,p)$ such that
 $$\left\|\int_0^te^{-(t-s)(-\Delta)^{\beta/2}} g(s)\,ds\right\|_{L^\infty(0,\infty;L^\phi(\mathbb{R}^n))}\leq C\|g\|_{L^\infty(0,\infty;L^1(\mathbb{R}^n)\cap L^{2p}(\mathbb{R}^n)\cap L^{\frac{2p}{p-1}}(\mathbb{R}^n))},$$ 
for every $g\in L^\infty(0,\infty;L^1(\mathbb{R}^n)\cap L^{2p}(\mathbb{R}^n)\cap L^{\frac{2p}{p-1}}(\mathbb{R}^n))$. 
\end{lemma}
\proof On the one hand, by (\ref{P_2}), we have
\begin{eqnarray*}
\int_{\mathbb{R}^n}\phi\left(\frac{|e^{-t(-\Delta)^{\beta/2}} \varphi|}{\lambda}\right)\,dx&=&\sum_{k=2}^{\infty}\frac{\left\|e^{-t(-\Delta)^{\beta/2}} \varphi\right\|_{pk}^{pk}}{k!\lambda^{pk}}\\
&\leq& \sum_{k=2}^{\infty}\frac{C^{pk}t^{-\frac{n}{\beta}(1-\frac{1}{pk})pk}\left\|\varphi\right\|_{1}^{pk}}{k!\lambda^{pk}}\\&=&t^{\frac{n}{\beta}}\phi\left(\frac{C t^{-\frac{n}{\beta}}\|\varphi\|_1}{\lambda}\right)\\
&\leq& t^{\frac{n}{\beta}}\left(\textnormal{exp}\left(\frac{C t^{-\frac{n}{\beta}}\|\varphi\|_1}{\lambda}\right)^{2p}-1\right),
\end{eqnarray*}
for all $t>0$, $\varphi\in L^1(\mathbb{R}^n)$, where we have used the fact that $e^{|x|^p}-1-|x|^p\leq e^{|x|^{2p}}-1$, for all $x\in\mathbb{R}$. As
$$ t^{\frac{n}{\beta}}\left(\textnormal{exp}\left(\frac{C t^{-\frac{n}{\beta}}\|\varphi\|_1}{\lambda}\right)^{2p}-1\right)\leq 1\Longleftrightarrow \lambda\geq C\,t^{-\frac{n}{\beta}}\left(\ln( t^{-\frac{n}{\beta}}+1)\right)^{-1/2p}\left\|\varphi\right\|_{1};$$
hence, 
\begin{eqnarray*}
&{}&\left\{\lambda>0,\,\, \lambda\in [ C\,t^{-\frac{n}{\beta}}\left(\ln( t^{-\frac{n}{\beta}}+1)\right)^{-1/2p}\left\|\varphi\right\|_{1};\infty[\right\}\\
&{}&\,\subseteq \left\{\lambda>0,\,\, \int_{\mathbb{R}^n}\phi\left(\frac{|e^{-t(-\Delta)^{\beta/2}} \varphi|}{\lambda}\right)\,dx\leq 1\right\};
\end{eqnarray*}
whereupon
\begin{equation}\label{eq12}
\left\|e^{-t(-\Delta)^{\beta/2}} \varphi\right\|_{L^\phi(\mathbb{R}^n)}\leq C\,t^{-\frac{n}{\beta}}\left(\ln( t^{-\frac{n}{\beta}}+1)\right)^{-1/2p}\left\|\varphi\right\|_{1},
\end{equation}
for all $t>0,$\, $\varphi\in L^1(\mathbb{R}^n)$. On the other hand, from (\ref{P_2}) and the embedding $L^{2p}(\mathbb{R}^n)\cap L^\infty(\mathbb{R}^n)\hookrightarrow L^\phi_0(\mathbb{R}^n)$ (see Lemma \ref{lemma2}), we have
\begin{eqnarray}\label{eq13}
\left\|e^{-t(-\Delta)^{\beta/2}} \varphi\right\|_{L^\phi(\mathbb{R}^n)}&\leq&\left\|e^{-t(-\Delta)^{\beta/2}} \varphi\right\|_{L^\infty(\mathbb{R}^n)}+\left\|e^{-t(-\Delta)^{\beta/2}} \varphi\right\|_{L^{2p}(\mathbb{R}^n)}\nonumber\\
&\leq& C t^{-\frac{n}{\beta}(\frac{p-1}{2p}-0)}\left\|\varphi\right\|_{L^{\frac{2p}{p-1}}(\mathbb{R}^n)}+\left\|\varphi\right\|_{L^{2p}(\mathbb{R}^n)}\nonumber\\
&=& C\,t^{-\frac{1}{2}}\left\|\varphi\right\|_{L^{\frac{2p}{p-1}}(\mathbb{R}^n)}+\left\|\varphi\right\|_{L^{2p}(\mathbb{R}^n)}\nonumber\\
&\leq& C\,(t^{-\frac{1}{2}}+1)\left(\left\|\varphi\right\|_{L^{\frac{2p}{p-1}}(\mathbb{R}^n)}+\left\|\varphi\right\|_{L^{2p}(\mathbb{R}^n)}\right),
\end{eqnarray}
for all $t>0$, $\varphi\in L^{2p}(\mathbb{R}^n)\cap L^{\frac{2p}{p-1}}(\mathbb{R}^n)$, where we have used the fact that $\beta=\frac{n(p-1)}{p}$. Now, let $g\in L^\infty(0,\infty;L^1(\mathbb{R}^n)\cap L^{2p}(\mathbb{R}^n)\cap L^{\frac{2p}{p-1}}(\mathbb{R}^n))$, we conclude from (\ref{eq12}) and (\ref{eq13}) that
$$
\left\|e^{-t(-\Delta)^{\beta/2}} g(t)\right\|_{L^\phi(\mathbb{R}^n)}\leq \kappa(t) \left\|g(t)\right\|_{L^1(\mathbb{R}^n)\cap L^{2p}(\mathbb{R}^n)\cap L^{\frac{2p}{p-1}}(\mathbb{R}^n)},
$$
for all $t>0$, where
$$\kappa(t) :=\min\left\{C\,(t^{-\frac{1}{2}}+1);C\,t^{-\frac{n}{\beta}}\left(\ln( t^{-\frac{n}{\beta}}+1)\right)^{-1/2p}\right\}.$$
We can easily check that $\kappa\in L^1(0,\infty)$. Therefore
\begin{eqnarray*}
\left\|\int_0^te^{-(t-s)(-\Delta)^{\beta/2}} g(s)\,ds\right\|_{L^\phi(\mathbb{R}^n)}&\leq& \int_0^t\left\|e^{-(t-s)(-\Delta)^{\beta/2}} g(s)\right\|_{L^\phi(\mathbb{R}^n)}\,ds\\
&\leq& \int_0^t\kappa(t-s)\left\|g(s)\right\|_{L^1(\mathbb{R}^n)\cap L^{2p}(\mathbb{R}^n)\cap L^{\frac{2p}{p-1}}(\mathbb{R}^n)}\,ds\\
&\leq&\|g\|_{L^\infty(0,\infty;L^1(\mathbb{R}^n)\cap L^{2p}(\mathbb{R}^n)\cap L^{\frac{2p}{p-1}}(\mathbb{R}^n))} \int_0^\infty\kappa(s)\,ds,
\end{eqnarray*}
for every $t>0$. This proves Lemma \ref{lemma5}. \hfill$\square$\\

Finally, the following proposition is needed for the local well-posedness result in the space $\textnormal{exp}\, L^p_0(\mathbb{R}^n)$.
\begin{proposition}\label{prop3}\cite[Proposition 2.9]{Majdoub2}\\
Let $1\leq p<\infty$ and $u\in C([0,T];\textnormal{exp}\, L^p_0(\mathbb{R}^n))$ for some $T>0$. Then, for every $\alpha>0$, it holds
$$\left(e^{\alpha|u|^p}-1\right)\in C([0,T];L^r(\mathbb{R}^n)),\quad 1\leq r<\infty.$$
\end{proposition}
\begin{corollary}\label{corollary1}\cite[Corollary 2.13]{Majdoub2}\\
Let $1\leq p<\infty$ and $u\in C([0,T];\textnormal{exp}\, L^p_0(\mathbb{R}^n))$ for some $T>0$. Assume that $f$ satisfies $(\ref{eq3})$. Then, for every $p\leq r<\infty$, it holds
$$f(u)\in C([0,T];L^r(\mathbb{R}^n)).$$
\end{corollary}


\section{Proof of Theorem \ref{theo1}}\label{sec3}
In this section, we prove Theorem \ref{theo1} i.e. the local existence and the uniqueness of a mild solution to (\ref{eq1}) in $C([0,T];\textnormal{exp}\, L^p_0(\mathbb{R}^n))$ for some $T>0$. Throughout this section, we assume that the nonlinearity $f$ satisfies (\ref{eq3}). In order to find the required solution, we will apply the Banach fixed-point theorem to the integral formulation (\ref{eq2}), using a decomposition argument developed in \cite{Ibrahim} and used in \cite{Ioku,Majdoub1,Majdoub2}. The idea is to split the initial data $u_0\in \textnormal{exp}\, L^p_0(\mathbb{R}^n)$, using the density of $C^\infty_0\mathbb{R}^n)$, into a small part in $\textnormal{exp}\, L^p(\mathbb{R}^n)$ and a smooth one. Let $u_0\in \textnormal{exp}\, L^p_0(\mathbb{R}^n)$. Then, for every $\varepsilon>0$ there exists $v_0\in C^\infty_0\mathbb{R}^n)$ such that
$$\|w_0\|_{\textnormal{exp}\, L^p(\mathbb{R}^n)}\leq\varepsilon,$$
where $w_0:=u_0-v_0$. Now, we split our problem (\ref{eq1}) into the following two problems. The first one is the fractional semilinear heat equation with smooth initial data:
\begin{equation}\label{eq24}
\left\{\begin{array}{ll}
\,\, \displaystyle {v_{t}+(-\Delta)^{\beta/2}v =f(v),} &\displaystyle {t>0,x\in {\mathbb{R}^n},}\\
{}\\
\displaystyle{v(0)=  v_0\in C^\infty_0(\mathbb{R}^n),\qquad\qquad}&\displaystyle{x\in {\mathbb{R}^n},}
\end{array}
\right.
\end{equation} 
and the second one is a fractional semilinear heat equation with small initial data in $\textnormal{exp}\, L^p(\mathbb{R}^n)$:
\begin{equation}\label{eq25}
\left\{\begin{array}{ll}
\,\, \displaystyle {w_{t}+(-\Delta)^{\beta/2}w =f(w+v)-f(v),} &\displaystyle {t>0,x\in {\mathbb{R}^n},}\\
{}\\
\displaystyle{w(0)= w_0,\,\,\|w_0\|_{\textnormal{exp}\, L^p}\leq\varepsilon,\qquad\qquad}&\displaystyle{x\in {\mathbb{R}^n}.}
\end{array}
\right.
\end{equation} 
We notice that if $v$ is a mild solution of (\ref{eq24}) and $w$ is a mild solution of (\ref{eq25}), then $u=v+w$ is a solution of our problem (\ref{eq2}), where the definition of the mild solutions for problems (\ref{eq24})- (\ref{eq25}) are defined similarly as in definition \ref{mild}. We now prove the local existence result concerning (\ref{eq24}) and (\ref{eq25}).
\begin{lemma}${}$\label{lemma7}\\
Let $0<\beta< 2$, $p>1$ and $v_0\in L^p(\mathbb{R}^n)\cap L^\infty(\mathbb{R}^n)$. Then, there exist a time $T=T(v_0)>0$ and a mild solution $v\in C([0,T];\textnormal{exp}\, L^p_0(\mathbb{R}^n))\cap L^\infty(0,T;L^\infty(\mathbb{R}^n))$ of $(\ref{eq24})$.
 \end{lemma}
\begin{lemma}${}$\label{lemma8}\\
Let $0<\beta< 2$, $p>1$, and $w_0\in \textnormal{exp}\, L^p_0(\mathbb{R}^n)$. Let $T>0$ and $v\in L^\infty(0,T;L^\infty(\mathbb{R}^n))$ be given by Lemma $\ref{lemma7}$. Then, for $\|w_0\|_{\textnormal{exp}\, L^p}\leq\varepsilon$, with $\varepsilon\ll1$ small enough, there exist a time $\widetilde{T}=\widetilde{T}(w_0,\varepsilon, v)>0$ and a mild solution $w\in C([0,\widetilde{T}];\textnormal{exp}\, L^p_0(\mathbb{R}^n))$ to problem $(\ref{eq25})$.
 \end{lemma}

\noindent {\bf Proof of Lemma \ref{lemma7}.} In order to use the Banach fixed-point theorem, we introduce the following complete metric space
$$Y_T:=\left\{v\in C([0,T];\textnormal{exp}\, L^p_0(\mathbb{R}^n))\cap L^\infty(0,T;L^\infty(\mathbb{R}^n));\,\,\|v\|_{Y_T}\leq 2\|v_0\|_{L^p\cap L^\infty}\right\},$$
where $\|v\|_{Y_T}:=\|v\|_{L^\infty(0,T;L^p)}+\|v\|_{L^\infty(0,T;L^\infty)}$ and $\|v_0\|_{L^p\cap L^\infty}:=\|v_0\|_{L^p}+\|v_0\|_{L^\infty}$. For $v\in Y_T$, we define $\Phi(v)$ by
$$\Phi(v):=e^{-t(-\Delta)^{\beta/2}}v_0+\int_0^te^{-(t-s)(-\Delta)^{\beta/2}}f(v(s))\,ds.$$
We will prove that if $T>0$ is small enough, then, $\Phi$ is a contraction from $Y_T$ into itself.\\
\noindent$\bullet$ {\bf $\Phi:Y_T\rightarrow Y_T$.} Let $v\in Y_T$. As $v_0\in L^p(\mathbb{R}^n)\cap L^\infty(\mathbb{R}^n)$, then, by Lemma \ref{lemma1}, we conclude that $v_0\in\textnormal{exp}\, L^p_0(\mathbb{R}^n)$. Then, using Proposition \ref{prop2}, $e^{-t(-\Delta)^{\beta/2}}v_0\in C([0,T];\textnormal{exp}\, L^p_0(\mathbb{R}^n))$. Next, for $q=p$ or $q=\infty$, we have
\begin{equation}\label{eq26}
\|f(v)\|_{L^q}\leq Ce^{\lambda\|v\|^p_\infty}\|v\|_{L^q}\leq Ce^{\lambda\|v\|^p_\infty}(2\|v_0\|_{L^p\cap L^\infty}),
\end{equation}
which implies, using again Lemma \ref{lemma1}, that $f(v)\in \textnormal{exp}\, L^p_0(\mathbb{R}^n)$ and more precisely $f(v)\in L^1(0,T;\textnormal{exp}\, L^p_0(\mathbb{R}^n))$ . It follows, by density and smoothing effect of the fractional semigroup $e^{-t(-\Delta)^{\beta/2}}$ (Lemma \ref{lemma10}), that 
$$\int_0^te^{-(t-s)(-\Delta)^{\beta/2}}f(v(s))\,ds\in C([0,T];\textnormal{exp}\, L^p_0(\mathbb{R}^n)).$$
So $\Phi(v)\in C([0,T];\textnormal{exp}\, L^p_0(\mathbb{R}^n))$. Moreover, using (\ref{P_3}) and (\ref{eq26}), we have
$$\|\Phi(v)\|_{Y_T}\leq \|v_0\|_{L^p\cap L^\infty} +2TC(2\|v_0\|_{L^p\cap L^\infty})e^{\lambda(2\|v_0\|_{L^p\cap L^\infty})^p}\leq 2\|v_0\|_{L^p\cap L^\infty},$$
for $T>0$ small enough, namely $4TCe^{\lambda(2\|v_0\|_{L^p\cap L^\infty})^p}\leq 1$. This proves that $\Phi(v)\in Y_T$.\\

\noindent$\bullet$ {\bf $\Phi$ is a contraction.} Let $v_1,v_2\in Y_T$. For $q=p$ or $q=\infty$, we have
\begin{eqnarray*}
\|f(v_1)-f(v_2)\|_{L^q}&\leq& C\|v_1-v_2\|_q(e^{\lambda\|v_1\|^p_\infty}+e^{\lambda\|v_2\|^p_\infty})\\
&\leq& 2C\|v_1-v_2\|_q e^{\lambda(2\|v_0\|_{L^p\cap L^\infty})^p}\\
&\leq& 2C\|v_1-v_2\|_{Y_T} e^{\lambda(2\|v_0\|_{L^p\cap L^\infty})^p}.
\end{eqnarray*}
By (\ref{P_3}), it holds
$$
\|\Phi(v_1)-\Phi(v_2)\|_{Y_T}\leq 2TC\|v_1-v_2\|_{Y_T} e^{\lambda(2\|v_0\|_{L^p\cap L^\infty})^p}\leq \frac{1}{2}\|v_1-v_2\|_{Y_T}.
$$
This finishes the proof of Lemma \ref{lemma7}. \hfill$\square$\\

\noindent {\bf Proof of Lemma \ref{lemma8}.} To prove Lemma \ref{lemma8}, we need the following result.
\begin{lemma}\cite[Lemma 4.4]{Majdoub2}\label{lemma9}\\
Let $v\in L^\infty(0,T;L^\infty(\mathbb{R}^n))$ for some $T>0$. Let $1<p\leq q<\infty$, and $w_1,w_2\in \textnormal{exp}\, L^p(\mathbb{R}^n)$ with $\|w_1\|_{\textnormal{exp}\, L^p},\|w_2\|_{\textnormal{exp}\, L^p}\leq M$ for sufficiently small $M>0$ (namely $2^p\lambda q M^p\leq1$, where $\lambda$ is given in (\ref{eq3})). Then, there exists a constant $C=C(q)>0$ such that 
$$\|f(w_1+v)-f(w_2+v)\|_q\leq Ce^{2^{p-1}\lambda\|v\|_\infty^p}\|w_1-w_2\|_{\textnormal{exp}\, L^p}.$$
 \end{lemma}
For $\widetilde{T}>0$, we define 
$$W_{\widetilde{T}}:=\left\{w\in C([0,\widetilde{T}];\textnormal{exp}\, L^p_0(\mathbb{R}^n));\,\,\|w\|_{L^\infty(0,\widetilde{T};\textnormal{exp}\, L^p_0)}\leq 2\varepsilon\right\},$$
and we consider the map $\Phi$ defined, for $w\in W_{\widetilde{T}}$, by
$$\Phi(w):=e^{-t(-\Delta)^{\beta/2}}w_0+\int_0^te^{-(t-s)(-\Delta)^{\beta/2}}\left(f(w(s)+v(s))-f(v(s))\right)\,ds.$$
We will prove that if $\varepsilon$ and $\widetilde{T}>0$ are small enough, then, $\Phi$ is a contraction from $W_{\widetilde{T}}$ into itself.\\

\noindent$\bullet$ {\bf $\Phi$ is a contraction.} Let $w_1,w_2\in W_{\widetilde{T}}$. Using Lemma \ref{lemma1}, i.e. the embedding $L^p(\mathbb{R}^n)\cap L^\infty(\mathbb{R}^n)\hookrightarrow \textnormal{exp}\, L^p_0(\mathbb{R}^n)$, we have
\begin{equation}\label{eq27}
\|\Phi(w_1)-\Phi(w_2)\|_{\textnormal{exp}\, L^p}\leq \frac{1}{(\ln 2)^{1/p}}\left(\|\Phi(w_1)-\Phi(w_2)\|_{p}+\|\Phi(w_1)-\Phi(w_2)\|_{\infty}\right).
\end{equation}
Let $r>0$ be an auxiliary constant such that $r>\max\{p,\frac{n}{\beta}\}$. Then
$$\|\Phi(w_1)-\Phi(w_2)\|_{\infty}\leq C\int_0^t(t-s)^{-\frac{n}{\beta r}}\|f(w_1(s)+v(s))-f(w_2(s)+v(s))\|_r\,ds,$$
thanks to the $L^r-L^\infty$ estimate (\ref{P_2}). Applying Lemma \ref{lemma9} with $q=r$ and under the condition $2^p\lambda r (2\varepsilon)^p\leq1$, we obtain
\begin{eqnarray}\label{eq28}
\|\Phi(w_1)-\Phi(w_2)\|_{\infty}&\leq& Ce^{2^{p-1}\lambda\|v\|_\infty^p}\left(\int_0^t(t-s)^{-\frac{n}{\beta r}}\,ds\right)\|w_1-w_2\|_{L^\infty(0,\widetilde{T};\textnormal{exp}\, L^p)}\nonumber\\
&\leq& Ce^{2^{p-1}\lambda\|v\|_\infty^p}\widetilde{T}^{1-\frac{n}{\beta r}}\|w_1-w_2\|_{L^\infty(0,\widetilde{T};\textnormal{exp}\, L^p)}.
\end{eqnarray}
On the other hand, applying again the $L^p-L^p$ estimate (\ref{P_3}), and Lemma \ref{lemma9} with $q=p$ under the condition $2^p\lambda p (2\varepsilon)^p\leq1$, we obtain
\begin{eqnarray}\label{eq29}
\|\Phi(w_1)-\Phi(w_2)\|_{p}&\leq& \int_0^t\left\|e^{-(t-s)(-\Delta)^{\beta/2}}\left(f(w_1(s)+v(s))-f(w_2(s)+v(s))\right)\right\|_{p}\,ds\nonumber\\
&\leq& \int_0^t\left\|f(w_1(s)+v(s))-f(w_2(s)+v(s))\right\|_{p}\,ds\nonumber\\
&\leq&Ce^{2^{p-1}\lambda\|v\|_\infty^p}\int_0^t\|w_1-w_2\|_{\textnormal{exp}\, L^p}\,ds\nonumber\\
&\leq& Ce^{2^{p-1}\lambda\|v\|_\infty^p}\widetilde{T}\|w_1-w_2\|_{L^\infty(0,\widetilde{T};\textnormal{exp}\, L^p)}.
\end{eqnarray}
Using (\ref{eq28}) and (\ref{eq29}) into (\ref{eq27}), we infer, by choosing $\varepsilon\ll1$ small enough, that
\begin{eqnarray}\label{eq30}
\|\Phi(w_1)-\Phi(w_2)\|_{\textnormal{exp}\, L^p}&\leq& Ce^{2^{p-1}\lambda\|v\|_\infty^p}\left(\widetilde{T}+\widetilde{T}^{1-\frac{n}{\beta r}}\right)\|w_1-w_2\|_{L^\infty(0,\widetilde{T};\textnormal{exp}\, L^p)}\nonumber\\
&\leq& \frac{1}{2}\|w_1-w_2\|_{L^\infty(0,\widetilde{T};\textnormal{exp}\, L^p)},
\end{eqnarray}
where $\widetilde{T}\ll1$ is chosen small enough such that $Ce^{2^{p-1}\lambda\|v\|_\infty^p}\left(\widetilde{T}+\widetilde{T}^{1-\frac{n}{\beta r}}\right)\leq\frac{1}{2}$.\\

\noindent$\bullet$ {\bf $\Phi:W_{\widetilde{T}}\rightarrow W_{\widetilde{T}}$}. Let $w\in W_{\widetilde{T}}$. As $w_0\in L^p(\mathbb{R}^n)\cap L^\infty(\mathbb{R}^n)$, then by Lemma \ref{lemma1}, we conclude that $w_0\in\textnormal{exp}\, L^p_0(\mathbb{R}^n)$. Then, using Proposition \ref{prop2}, 
$$e^{-t(-\Delta)^{\beta/2}}w_0\in C([0,T];\textnormal{exp}\, L^p_0(\mathbb{R}^n)).$$ 
Next, the estimates (\ref{eq28})-(\ref{eq29}) with $w_1=w$ and $w_2=0$, under the condition $2^p\lambda r (2\varepsilon)^p\leq1$, show that the nonlinear term satisfies
$$\Phi(w)-e^{-t(-\Delta)^{\beta/2}}w_0\in L^\infty(0,T;\textnormal{exp}\, L^p_0(\mathbb{R}^n)),$$
thanks to the embedding $L^p(\mathbb{R}^n)\cap L^\infty(\mathbb{R}^n)\hookrightarrow \textnormal{exp}\, L^p_0(\mathbb{R}^n)$ (Lemma \ref{lemma1}). By the standard smoothing effect of the fractional semigroup $e^{-t(-\Delta)^{\beta/2}}$ (Lemma \ref{lemma10}), it follows that 
$$\Phi(w)-e^{-t(-\Delta)^{\beta/2}}w_0\in C([0,T];\textnormal{exp}\, L^p_0(\mathbb{R}^n)).$$
So 
$$\Phi(w)\in C([0,T];\textnormal{exp}\, L^p_0(\mathbb{R}^n)).$$ 
Moreover, using Proposition \ref{prop1}, and (\ref{eq30}) with $w_1=w$ and $w_2=0$ for $T\ll1$, we have
$$\|\Phi(w)\|_{W_{\widetilde{T}}}\leq \|w_0\|_{\textnormal{exp}\, L^p} +\frac{1}{2}\|w\|_{L^\infty(0,\widetilde{T};\textnormal{exp}\, L^p)}\leq \varepsilon+\frac{1}{2}(2\varepsilon)=2\varepsilon.$$
 This proves that $\Phi(w)\in W_{\widetilde{T}}$.\hfill$\square$\\

\noindent {\bf Proof of the existence part in Theorem \ref{theo1}.} We choose $T$, $\varepsilon$, and $\widetilde{T}$ in the following order. Let $r>\max\{p,\frac{n}{\beta}\}$ and fix $\varepsilon>0$ such that 
$$2^p\lambda r (2\varepsilon)^p\leq1.$$
Next, one can decompose $u_0=v_0+w_0$ with $v_0\in C^\infty_0\mathbb{R}^n)$ and $\|w_0\|_{\textnormal{exp}\, L^p(\mathbb{R}^n)}\leq\varepsilon$. By Lemma \ref{lemma7}, there exist a time $0<T_1=T_1(\|v_0\|_{L^p\cap L^\infty})\ll1$ and a mild solution $v\in C([0,T_1];\textnormal{exp}\, L^p_0(\mathbb{R}^n))\cap L^\infty(0,T_1;L^\infty(\mathbb{R}^n))$ of (\ref{eq24}) such that  $\|v\|_{L^\infty(0,T;L^p\cap L^\infty)}\leq 2\|v_0\|_{L^p\cap L^\infty}$. By Choosing  $\widetilde{T}>0$ small enough such that $\widetilde{T}<T_1$ and
$$Ce^{2^{2p-1}\lambda\|v_0\|_{L^p\cap L^\infty}^p}\left(\widetilde{T}+\widetilde{T}^{1-\frac{n}{\beta r}}\right)\leq\frac{1}{2},$$
and using Lemma \ref{lemma8}, there exists a mild solution $w\in C([0,\widetilde{T}];\textnormal{exp}\, L^p_0(\mathbb{R}^n))$ to problem (\ref{eq25}). We conclude that $u:v+w$ is a mild solution of (\ref{eq1}) in 
$C([0,\widetilde{T}];\textnormal{exp}\, L^p_0(\mathbb{R}^n))$.\hfill$\square$\\

\noindent {\bf Proof of the uniqueness part in Theorem \ref{theo1}.} Let us suppose that $u,v\in C([0,T];\textnormal{exp}\, L^p_0(\mathbb{R}^n))$ are two mild solutions of (\ref{eq1}) for some $T>0$, and with the same initial data $u(0)=v(0)=u_0$. Let
$$t_0=\sup\{t\in[0,T]\quad \mbox{such that}\quad u(s)=v(s)\quad\mbox{for every}\,\,s\in[0,t]\}.$$
Let us suppose  that $0\leq t_0<T$. Since $u(t)$ and $v(t)$ are continuous in time, we have $u(t_0)=v(t_0)$. Let us denote $\tilde{u}(t):=u(t+t_0)$ and $\tilde{v}(t):=v(t+t_0)$. Then $\tilde{u},\tilde{v}\in C([0,T-t_0];\textnormal{exp}\, L^p_0(\mathbb{R}^n))$ and satisfy (\ref{eq2}) on $(0,T-t_0]$ with $\tilde{u}(0)=\tilde{v}(0)=u(t_0)$. We will prove that there exists a short positive time $0<\tilde{t}\leq T-t_0$ such that 
\begin{equation}\label{eq31}
\sup_{0<t<\tilde{t}}\|\tilde{u}(t)-\tilde{v}(t)\|_{\textnormal{exp}\, L^p}\leq C(\tilde{t})\sup_{0<t<\tilde{t}}\|\tilde{u}(t)-\tilde{v}(t)\|_{\textnormal{exp}\, L^p},
\end{equation}
for some $C(\tilde{t})<1$, and so $\tilde{u}(t)=\tilde{v}(t)$ for any $t\in[0,\tilde{t}]$. Therefore $u(t+t_0)=v(t+t_0)$ for any $t\in[0,\tilde{t}]$ which is a contradiction with the definition of $t_0$. In order to establish inequality (\ref{eq31}), we control both the $L^p$-norm and $L^\infty$-norm of $\tilde{u}-\tilde{v}$. Using $L^p-L^p$ estimate (\ref{P_3}), we obtain
\begin{eqnarray*}
\|\tilde{u}(t)-\tilde{v}(t)\|_p&\leq& \int_0^t\left\|e^{-(t-s)(-\Delta)^{\beta/2}}(f(\tilde{u}(s))-f(\tilde{v}(s)))\right\|_p\,ds\\
&\leq& \int_0^t\|(f(\tilde{u}(s))-f(\tilde{v}(s)))\|_p\,ds.
\end{eqnarray*}
By (\ref{eq3}) and H\"older's inequality, we get
\begin{eqnarray*}
&{}&\|\tilde{u}(t)-\tilde{v}(t)\|_p\\
&{}&\,\leq C \int_0^t\|(\tilde{u}(s)-\tilde{v}(s))(e^{\lambda |\tilde{u}|^p}+e^{\lambda |\tilde{v}|^p})\|_p\,ds\\
&{}&\,\leq 2\int_0^t\|\tilde{u}(s)-\tilde{v}(s)\|_p\,ds+\int_0^t\left\|(\tilde{u}(s)-\tilde{v}(s))((e^{\lambda |\tilde{u}|^p}-1)+(e^{\lambda |\tilde{v}|^p}-1))\right\|_p\,ds\\
&{}&\,\leq 2\int_0^t\|\tilde{u}(s)-\tilde{v}(s)\|_p\,ds+\int_0^t\|\tilde{u}(s)-\tilde{v}(s)\|_q\left\|(e^{\lambda |\tilde{u}|^p}-1)+(e^{\lambda |\tilde{v}|^p}-1)\right\|_r\,ds,
\end{eqnarray*}
where $\frac{1}{q}+\frac{1}{r}=\frac{1}{p}$.
Thanks to Lemma \ref{lemma3} and $q\geq p$, we infer that
\begin{eqnarray*}
\|\tilde{u}(t)-\tilde{v}(t)\|_p&\leq& Ct\sup_{0<s<t}\|\tilde{u}(s)-\tilde{v}(s)\|_{\textnormal{exp}\, L^p}\\
&{}&+\,C\sup_{0<s<t}\|\tilde{u}(s)-\tilde{v}(s)\|_{\textnormal{exp}\, L^p}\int_0^t\left\|(e^{\lambda |\tilde{u}|^p}-1)+(e^{\lambda |\tilde{v}|^p}-1)\right\|_r\,ds.
\end{eqnarray*}
Moreover, using Proposition \ref{prop3}, we obtain
\begin{eqnarray}\label{eq33}
&{}&\,\sup_{0<s<T-t_0}\left\|(e^{\lambda |\tilde{u}|^p}-1)+(e^{\lambda |\tilde{v}|^p}-1)\right\|_r\nonumber\\
&{}&\,\leq \sup_{0<s<T-t_0}\left[\left(\int_{\mathbb{R}^n}(e^{\lambda r |\tilde{u}|^p}-1)\,dx\right)^{1/r}+\left(\int_{\mathbb{R}^n}(e^{\lambda r |\tilde{v}|^p}-1)\,dx\right)^{1/r}\right]\nonumber\\
&{}&\,\leq C(T,t_0,\tilde{u},\tilde{v})<\infty.
\end{eqnarray}
Consequently,
\begin{equation}\label{eq32}
\sup_{0<s<t}\|\tilde{u}(s)-\tilde{v}(s)\|_{ L^p}\leq C(T,t_0,\tilde{u},\tilde{v})t \sup_{0<s<t}\|\tilde{u}(s)-\tilde{v}(s)\|_{\textnormal{exp}\, L^p}.
\end{equation}
In a similar way, using $L^r-L^\infty$ estimate (\ref{P_2}), we obtain
\begin{eqnarray*}
\|\tilde{u}(t)-\tilde{v}(t)\|_\infty&\leq& \int_0^t\left\|e^{-(t-s)(-\Delta)^{\beta/2}}(f(\tilde{u}(s))-f(\tilde{v}(s)))\right\|_\infty\,ds\\
&\leq& C\int_0^t(t-s)^{-\frac{n}{\beta r}}\|(f(\tilde{u}(s))-f(\tilde{v}(s)))\|_r\,ds,
\end{eqnarray*}
for some $r>\max\{p,\frac{n}{\beta}\}$. By (\ref{eq3}) and H\"older's inequality, we get
\begin{eqnarray*}
\|\tilde{u}(t)-\tilde{v}(t)\|_\infty&\leq& C \int_0^t(t-s)^{-\frac{n}{\beta r}}\|(\tilde{u}(s)-\tilde{v}(s))(e^{\lambda |\tilde{u}|^p}+e^{\lambda |\tilde{v}|^p})\|_r\,ds\\
&\leq& 2\int_0^t(t-s)^{-\frac{n}{\beta r}}\|\tilde{u}(s)-\tilde{v}(s)\|_r\,ds\\
&{}&\,+\int_0^t(t-s)^{-\frac{n}{\beta r}}\left\|(\tilde{u}(s)-\tilde{v}(s))((e^{\lambda |\tilde{u}|^p}-1)+(e^{\lambda |\tilde{v}|^p}-1))\right\|_r\,ds\\
&\leq& 2\int_0^t(t-s)^{-\frac{n}{\beta r}}\|\tilde{u}(s)-\tilde{v}(s)\|_r\,ds\\
&{}&\,+\int_0^t(t-s)^{-\frac{n}{\beta r}}\|\tilde{u}(s)-\tilde{v}(s)\|_{\tilde{q}}\left\|(e^{\lambda |\tilde{u}|^p}-1)+(e^{\lambda |\tilde{v}|^p}-1)\right\|_{\tilde{r}}\,ds,
\end{eqnarray*}
where $\frac{1}{\tilde{q}}+\frac{1}{\tilde{r}}=\frac{1}{r}$. Since $\tilde{q},\tilde{r}\geq r> p$, one can apply an estimate similar to (\ref{eq33}) via Lemma \ref{lemma3} and Proposition \ref{prop3}, and obtain that
\begin{equation}\label{eq34}
\sup_{0<s<t}\|\tilde{u}(s)-\tilde{v}(s)\|_{ L^\infty}\leq C(T,t_0,\tilde{u},\tilde{v})t^{1-\frac{n}{\beta r}} \sup_{0<s<t}\|\tilde{u}(s)-\tilde{v}(s)\|_{\textnormal{exp}\, L^p}.
\end{equation}
Finally, the two inequalities (\ref{eq32}) and (\ref{eq34}) with the embedding relation $L^p(\mathbb{R}^n)\cap L^\infty(\mathbb{R}^n)\hookrightarrow \textnormal{exp}\, L^p_0(\mathbb{R}^n)$ (Lemma \ref{lemma1}) imply
$$\sup_{0<s<t}\|\tilde{u}(s)-\tilde{v}(s)\|_{\textnormal{exp}\, L^p}\leq C(T,t_0,\tilde{u},\tilde{v})(t+t^{1-\frac{n}{\beta r}}) \sup_{0<s<t}\|\tilde{u}(s)-\tilde{v}(s)\|_{\textnormal{exp}\, L^p},$$
and for $t$ small enough, we obtain the desired estimate (\ref{eq31}). \hfill$\square$\\

\begin{rmk}
The solution in Theorem \ref{theo1} belongs to $L^\infty_{loc}(0,T;L^\infty(\mathbb{R}^n))$. Indeed, let $u\in C([0,T];\textnormal{exp}\, L^p_0(\mathbb{R}^n))$ be a mild solution of (\ref{eq1}) i.e. a solution of the integral equation (\ref{eq2}). Using $L^p-L^\infty$ estimate (\ref{P_2}) and Lemma \ref{lemma3}, we get 
$$\|e^{-t(-\Delta)^{\beta/2}}u_0\|_\infty\leq Ct^{-\frac{n}{\beta p}}\|u_0\|_p\leq Ct^{-\frac{n}{\beta p}}\|u_0\|_{\textnormal{exp}\, L^p},$$
for all $0<t<T$. Hence $e^{-t(-\Delta)^{\beta/2}}u_0\in L^\infty(\mathbb{R}^n)$ for all $0<t<T$. Thus it remains to estimate the nonlinear term. Fix $r>\max\{p,\frac{n}{\beta}\}$, using $L^r-L^\infty$ estimate (\ref{P_2}), we get
\begin{eqnarray*}
\left\|\int_0^te^{-(t-s)(-\Delta)^{\beta/2}}f(u(s))\,ds\right\|_\infty&\leq&\int_0^t \left\|e^{-(t-s)(-\Delta)^{\beta/2}}f(u(s))\right\|_\infty\,ds\\
&\leq& \int_0^t t-s)^{-\frac{n}{\beta r}}\|f(u(s))\|_r\,ds\\
&\leq& Ct^{1-\frac{n}{\beta r}}\sup_{0\leq t\leq T}\|f(u(t))\|_r<\infty,
\end{eqnarray*}
where we have used Corollary \ref{corollary1}. This shows that $u\in L^\infty_{loc}(0,T;L^\infty(\mathbb{R}^n))$. In particular, if $f\in C^1(\mathbb{R}^n)$, the solution $u\in  C([0,T];\textnormal{exp}\, L^p_0(\mathbb{R}^n))\cap L^\infty_{loc}(0,T;L^\infty(\mathbb{R}^n))$ satisfies (\ref{eq1}) in the classical sense, i.e. $C^1$ in time $t\in(0,T)$ and $C^2$ in space $\mathbb{R}^n$.
\end{rmk}

\begin{rmk}
Using the uniqueness, the constructed solution $u$ of (\ref{eq1}) can be extended to a maximal interval $[0,T_{\max})$ by well known argument (see cf. Cazenave et Haraux \cite{CH}) where
$$
T_{\max}:=\sup\left\{T>0\;;\;\textnormal{there exist a mild solution} \,\, u\in C([0,T];\textnormal{exp}\, L^p_0(\mathbb{R}^n))\,\hbox{to}\,\, (\ref{eq1})\right\}\leq+\infty.
$$
Moreover, if $T_{\max}<\infty$, then 
$$\lim_{t\rightarrow T_{\max}}\|u(t)\|_{L^p\cap L^\infty(\mathbb{R}^n)}=\infty.$$
\end{rmk}


\section{Proof of Theorem \ref{theo2}}\label{sec4}

\subsection{\bf {Proof of global existence in Theorem \ref{theo2} (case of $\beta<\frac{n(p-1)}{p}$)}}\label{subsec4.1}
In this subsection, we prove the global existence of solution in Theorem \ref{theo2} in the case of $\beta<\frac{n(p-1)}{p}$. We will use the fixed-point theorem. Let us first define the following space
$$E_\varepsilon=\left\{u\in L^\infty(0,\infty;\textnormal{exp~$L^p(\mathbb{R}^n)$});\,\,\|u\|_{L^\infty(0,\infty;\textnormal{exp~$L^p(\mathbb{R}^n)$})}\leq 2\varepsilon\right\},$$
where $\varepsilon>0$ is a positive constant, small enough, that will be chosen later such that $\|u_0\|_{\textnormal{exp~$L^p(\mathbb{R}^n)$}}\leq\varepsilon$. For $u\in E_\varepsilon$, we define $\Phi(u)$ by
$$\Phi(u):=e^{-t(-\Delta)^{\beta/2}}u_0+\int_0^te^{-(t-s)(-\Delta)^{\beta/2}}f(u(s))\,ds.$$
Our goal is to prove that $\Phi:E_\varepsilon\rightarrow E_\varepsilon$ is a contraction map.\\

\noindent$\bullet$ {\bf $\Phi:E_\varepsilon\rightarrow E_\varepsilon$.} Let $u\in E_\varepsilon$, we have
\begin{eqnarray*}
\left\|\Phi(u)\right\|_{\textnormal{exp~$L^p(\mathbb{R}^n)$}}&\leq&\left\|e^{-t(-\Delta)^{\beta/2}}u_0\right\|_{\textnormal{exp~$L^p(\mathbb{R}^n)$}}+\left\|\int_0^te^{-(t-s)(-\Delta)^{\beta/2}} f(u(s))\,ds\right\|_{\textnormal{exp~$L^p(\mathbb{R}^n)$}}\\
&\leq&\left\|u_0\right\|_{\textnormal{exp~$L^p(\mathbb{R}^n)$}}+C\left\|f(u)\right\|_{L^\infty(0,\infty;L^1(\mathbb{R}^n)\cap L^r(\mathbb{R}^n))}\\
&\leq&\varepsilon+C\left\|f(u)\right\|_{L^\infty(0,\infty;L^1(\mathbb{R}^n)\cap L^r(\mathbb{R}^n))},
\end{eqnarray*}
for every $r>\frac{n}{\beta}>1$, where we have used Proposition \ref{prop1} and Lemma \ref{lemma4}. It remains to estimate $f(u)$ in $L^q(\mathbb{R}^n)$ for $q=1,r$. From the assumption (\ref{eq4}), we see
$$|f(u)|\leq C|u|^{m}e^{\lambda |u|^p}= C|u|^{m}\left(e^{\lambda |u|^p}-1\right)+C|u|^{m},$$
then, by H\"older's inequality, we obtain
\begin{eqnarray*}
\left\|f(u)\right\|_{L^q(\mathbb{R}^n)}&\leq& C \left\|u\right\|^{m}_{L^{2mq}(\mathbb{R}^n)} \left\|e^{\lambda |u|^p}-1\right\|_{L^{2q}(\mathbb{R}^n)}+C\left\|u\right\|^{m}_{L^{mq}(\mathbb{R}^n)}\\
&\leq& C \left\|u\right\|^{m}_{\textnormal{exp~$L^p(\mathbb{R}^n)$}} \left\|e^{\lambda |u|^p}-1\right\|_{L^{2q}(\mathbb{R}^n)}+C\left\|u\right\|^{m}_{\textnormal{exp~$L^p(\mathbb{R}^n)$}},
\end{eqnarray*}
where we have used Lemma \ref{lemma3} and $m\geq p$. Next, using Lemma \ref{lemma6} and the fact that $u\in E_\varepsilon$, we have
\begin{equation}\label{eq19}
\left\|f(u)\right\|_{L^q(\mathbb{R}^n)}\leq C(2\varepsilon)^m\left(1+(2\lambda q(2\varepsilon)^p)^{1/q}\right)\leq C(2\varepsilon)^m\left(1+(2\lambda q(2\varepsilon)^p)^{1/r}\right).
\end{equation}
This implies, by choosing $\varepsilon$ small enough, that
$$\left\|\Phi(u)\right\|_{\textnormal{exp~$L^p(\mathbb{R}^n)$}}\leq \varepsilon+C(2\varepsilon)^m\left(1+(2\lambda q(2\varepsilon)^p)^{1/r}\right)\leq 2\varepsilon,$$
i.e. $\Phi(u)\in E_\varepsilon$.\\

\noindent$\bullet$ {\bf $\Phi$ is a contraction.} Let $u,v\in E_\varepsilon$, we have
\begin{eqnarray*}
\left\|\Phi(u)-\Phi(v)\right\|_{\textnormal{exp~$L^p(\mathbb{R}^n)$}}&=&\left\|\int_0^te^{-(t-s)(-\Delta)^{\beta/2}} \left(f(u(s))-f(v(s))\right)\,ds\right\|_{\textnormal{exp~$L^p(\mathbb{R}^n)$}}\\
&\leq&C\left\|f(u)-f(v)\right\|_{L^\infty(0,\infty;L^1(\mathbb{R}^n)\cap L^r(\mathbb{R}^n))},
\end{eqnarray*}
for every $r>\frac{n}{\beta}>1$, where we have used Lemma \ref{lemma4}. To estimate $f(u)-f(v)$ in $L^1(\mathbb{R}^n)\cap L^r(\mathbb{R}^n)$, let $q=1,r$. We see, using assumption (\ref{eq4}), that
\begin{eqnarray*}
|f(u)-f(v)|&\leq& C|u-v|\left(|u|^{m-1}e^{\lambda |u|^p}+|v|^{m-1}e^{\lambda |v|^p}\right)\\
&=& C|u-v|\left(|u|^{m-1}\left(e^{\lambda |u|^p}-1\right)+|v|^{m-1}\left(e^{\lambda |v|^p}-1\right)\right)\\
&{}&+\,C|u-v|\left(|u|^{m-1}+|v|^{m-1}\right),
\end{eqnarray*}
then, by H\"older's inequality, we obtain
$$C\left\|f(u)-f(v)\right\|_{L^q(\mathbb{R}^n)}\leq I+II,$$
where
$$I:=C \left\|u-v\right\|_{L^{mq}(\mathbb{R}^n)} \left\||u|^{m-1}\left(e^{\lambda |u|^p}-1\right)+|v|^{m-1}\left(e^{\lambda |v|^p}-1\right)\right\|_{L^{\frac{mq}{m-1}}(\mathbb{R}^n)},$$
and
$$
II:=C\left\|u-v\right\|_{L^{mq}(\mathbb{R}^n)}\left\||u|^{m-1}+|v|^{m-1}\right\|_{L^{\frac{mq}{m-1}}(\mathbb{R}^n)}.
$$
Using again H\"older's inequality, Lemma \ref{lemma3}, and  $m\geq p$, we get
\begin{eqnarray*}
I&\leq& C \left\|u-v\right\|_{\textnormal{exp~$L^p(\mathbb{R}^n)$}}\\
&{}&\times \left( \left\||u|^{m-1}\left(e^{\lambda |u|^p}-1\right)\right\|_{L^{\frac{mq}{m-1}}(\mathbb{R}^n)}+\left\||v|^{m-1}\left(e^{\lambda |v|^p}-1\right)\right\|_{L^{\frac{mq}{m-1}}(\mathbb{R}^n)}\right)\\
&\leq&C \left\|u-v\right\|_{\textnormal{exp~$L^p(\mathbb{R}^n)$}}\\
&{}&\times\left( \|u\|^{m-1}_{L^{2mq}(\mathbb{R}^n)}\left\|e^{\lambda |u|^p}-1\right\|_{L^{\frac{2mq}{m-1}}(\mathbb{R}^n)}+\|v\|^{m-1}_{L^{2mq}(\mathbb{R}^n)}\left\|e^{\lambda |v|^p}-1\right\|_{L^{\frac{2mq}{m-1}}(\mathbb{R}^n)}\right)\\
&\leq&C \left\|u-v\right\|_{\textnormal{exp~$L^p(\mathbb{R}^n)$}}\\
&{}&\times\left( \|u\|^{m-1}_{\textnormal{exp~$L^p(\mathbb{R}^n)$}}\left\|e^{\lambda |u|^p}-1\right\|_{L^{\frac{2mq}{m-1}}(\mathbb{R}^n)}+\|v\|^{m-1}_{\textnormal{exp~$L^p(\mathbb{R}^n)$}}\left\|e^{\lambda |v|^p}-1\right\|_{L^{\frac{2mq}{m-1}}(\mathbb{R}^n)}\right).
\end{eqnarray*}
Then, using Lemma \ref{lemma6} and the fact that $u,v\in E_\varepsilon$, we have
$$I\leq C2^m\varepsilon^{m-1}\left(\frac{2\lambda mq}{m-1}(2\varepsilon)^p\right)^{\frac{m-1}{2mq}}\left\|u-v\right\|_{E_\varepsilon}\leq \frac{1}{8}\left\|u-v\right\|_{E_\varepsilon},$$
for $\varepsilon>0$ small enough. Similarly,
\begin{eqnarray*}
II&\leq& C\left\|u-v\right\|_{\textnormal{exp~$L^p(\mathbb{R}^n)$}}\left(\|u\|_{L^{mq}(\mathbb{R}^n)}^{m-1}+\|v\|_{L^{mq}(\mathbb{R}^n)}^{m-1}\right)\\
&\leq& C\left\|u-v\right\|_{\textnormal{exp~$L^p(\mathbb{R}^n)$}}\left(\|u\|_{\textnormal{exp~$L^p(\mathbb{R}^n)$}}^{m-1}+\|v\|_{\textnormal{exp~$L^p(\mathbb{R}^n)$}}^{m-1}\right)\\
&\leq& C2^m\varepsilon^{m-1}\left\|u-v\right\|_{\textnormal{exp~$L^p(\mathbb{R}^n)$}}\\
&\leq&  \frac{1}{8}\left\|u-v\right\|_{E_\varepsilon},
\end{eqnarray*}
for $\varepsilon>0$ small enough. We conclude that
$$C\left\|f(u)-f(v)\right\|_{L^1(\mathbb{R}^n)\cap L^r(\mathbb{R}^n)}\leq 2(I+II).$$
Hence,
$$\left\|\Phi(u)-\Phi(v)\right\|_{\textnormal{exp~$L^p(\mathbb{R}^n)$}}\leq  \frac{1}{2}\left\|u-v\right\|_{E_\varepsilon}.$$
This completes the proof of the existence of global solution in Theorem \ref{theo2} in the case of $\beta<\frac{n(p-1)}{p}$. To obtain the decay estimate (\ref{eq35}), we follow the same calculation as in the part of contraction mapping in the Subsection \ref{subsec4.2} below where we consider, instead of the space $E_\varepsilon$, the following complete metric space
$$\left\{u\in L^\infty(0,\infty;\textnormal{exp~$L^p(\mathbb{R}^n)$});\,\,\sup_{t>0}t^\sigma\|u(t)\|_{L^q(\mathbb{R}^n)}+\|u\|_{L^\infty(0,\infty;\textnormal{exp~$L^p(\mathbb{R}^n)$})}\leq M\varepsilon\right\},$$
endowed by the metric $d$ defined by $d(u,v):=\sup_{t>0}t^\sigma\|u(t)-v(t)\|_{L^q(\mathbb{R}^n)}$,
for certain large constant $M>0$, where $0<\varepsilon\ll1$ is a positive constant, small enough, that will be chosen later such that $\|u_0\|_{\textnormal{exp~$L^p(\mathbb{R}^n)$}}\leq\varepsilon$. The new parameters $\sigma$ and $q$ are chosen as follows:
$$ \sigma=\frac{1}{m-1}-\frac{n}{\beta q}>0,$$
and
$$ \frac{n(m-1)}{\beta}<q<\frac{n(m-1)}{\beta}\frac{1}{(2-m)_+}.$$
\hfill$\square$\\


\subsection{\bf {Proof of global existence in Theorem \ref{theo2} (case of $\beta\geq\frac{n(p-1)}{p}$)}}\label{subsec4.2}
This subsection is devoted to prove the existence of global solution in Theorem \ref{theo2} in the case of $\beta\geq\frac{n(p-1)}{p}$ by using same ideas as in \cite{Majdoub2} together with Lemma \ref{lemma5}. As the last section, we will use a contraction mapping argument in an appropriate complete space. Let us define
$$B_\varepsilon=\left\{u\in L^\infty(0,\infty;\textnormal{exp~$L^p(\mathbb{R}^n)$});\,\,\sup_{t>0}t^\sigma\|u(t)\|_{L^q(\mathbb{R}^n)}+\|u\|_{L^\infty(0,\infty;\textnormal{exp~$L^p(\mathbb{R}^n)$})}\leq M\varepsilon\right\},$$
for certain large constant $M>0$, where $0<\varepsilon\ll1$ is a positive constant, small enough, that will be chosen later such that $\|u_0\|_{\textnormal{exp~$L^p(\mathbb{R}^n)$}}\leq\varepsilon$. Using Proposition 2.2 in \cite{Majdoub2}, we can check that $B_\varepsilon$ is a complete metric space with the distance $d(u,v):=\sup_{t>0}t^\sigma\|u(t)-v(t)\|_{L^q(\mathbb{R}^n)}$. For $u\in B_\varepsilon$, we define, as above, $\Phi(u)$ by
$$\Phi(u):=e^{-t(-\Delta)^{\beta/2}}u_0+\int_0^te^{-(t-s)(-\Delta)^{\beta/2}}f(u(s))\,ds.$$

\noindent$\bullet$ {\bf $\Phi:B_\varepsilon\rightarrow B_\varepsilon$ .} Let $u\in B_\varepsilon$. By Proposition \ref{prop1}, we have
$$\|e^{-t(-\Delta)^{\beta/2}}u_0\|_{\textnormal{exp~$L^p(\mathbb{R}^n)$}}\leq \|u_0\|_{\textnormal{exp~$L^p(\mathbb{R}^n)$}}\leq\varepsilon.$$
Moreover, by choosing $\sigma=\frac{n}{\beta}\left(\frac{\beta}{n(m-1)}-\frac{1}{q}\right)=\frac{1}{m-1}-\frac{n}{\beta q}>0$, for $q>\frac{n(m-1)}{\beta}\geq p$, and using Lemma \ref{lemma3}, we get
\begin{eqnarray*}
t^\sigma\|e^{-t(-\Delta)^{\beta/2}}u_0\|_{L^q(\mathbb{R}^n)}&\leq& C t^\sigma t^{-\frac{n}{\beta}\left(\frac{\beta}{n(m-1)}-\frac{1}{q}\right)}\|u_0\|_{L^{\frac{n(m-1)}{\beta}}(\mathbb{R}^n)}\\
&=&C\,\|u_0\|_{L^{\frac{n(m-1)}{\beta}}(\mathbb{R}^n)}\\
&\leq& C\,\|u_0\|_{\textnormal{exp~$L^p(\mathbb{R}^n)$}}\\
&\leq& C\varepsilon.
\end{eqnarray*}
To estimate the second term in $\Phi(u)$ on $\textnormal{exp~$L^p(\mathbb{R}^n)$}$, we start to study the case of $\beta=\frac{n(p-1)}{p}$ by remembering (see (\ref{eq0})) that
$$C_1\|u\|_{\textnormal{exp~$L^p(\mathbb{R}^n)$}}\leq \|u\|_{L^p(\mathbb{R}^n)}+\|u\|_{L^\phi(\mathbb{R}^n)}  \leq C_2 \|u\|_{\textnormal{exp~$L^p(\mathbb{R}^n)$}},$$
for some $C_1,C_2>0$, where $\phi(u)=e^{|u|^p}-1-|u|^p$. Therefore, it is enough to prove the two following inequalities:
\begin{equation}\label{eq113}
\left\|\int_0^te^{-(t-s)(-\Delta)^{\beta/2}} f(u(s))\,ds\right\|_{L^\infty(0,\infty;L^p(\mathbb{R}^n))}=O(\varepsilon),
\end{equation}
and
\begin{equation}\label{eq14}
\left\|\int_0^te^{-(t-s)(-\Delta)^{\beta/2}}f(u(s))\,ds\right\|_{L^\infty(0,\infty;L^\phi(\mathbb{R}^n))}=O(\varepsilon).
\end{equation}
We start to prove (\ref{eq113}). As 
$$|f(u)|\leq C|u|^me^{\lambda|u|^p}=C|u|^m\sum_{k=0}^{\infty}\frac{\lambda^k}{k!}|u|^{kp}=C|u|\sum_{k=0}^{\infty}\frac{\lambda^k}{k!}|u|^{kp+m-1},$$
we have
\begin{eqnarray*}
&{}&\left\|\int_0^te^{-(t-s)(-\Delta)^{\beta/2}} f(u(s))\,ds\right\|_{L^p(\mathbb{R}^n)}\\
&{}&\leq C\int_0^t(t-s)^{-\frac{n}{\beta}\left(\frac{1}{r}-\frac{1}{p}\right)}\|f(u(s))\|_{L^r(\mathbb{R}^n)}\,ds\\
&{}&\leq C\sum_{k=0}^{\infty}\frac{\lambda^k}{k!}\int_0^t(t-s)^{-\frac{n}{\beta}\left(\frac{1}{r}-\frac{1}{p}\right)}\|u(s)\|_{L^p(\mathbb{R}^n)}\||u(s)|^{kp+m-1}\|_{L^a(\mathbb{R}^n)}\,ds\\
&{}&=C\sum_{k=0}^{\infty}\frac{\lambda^k}{k!}\int_0^t(t-s)^{-\frac{n}{\beta}\left(\frac{1}{r}-\frac{1}{p}\right)}\|u(s)\|_{L^p(\mathbb{R}^n)}\|u(s)\|^{kp+m-1}_{L^{(kp+m-1)a}(\mathbb{R}^n)}\,ds,
\end{eqnarray*}
where we have used (\ref{P_2}) and H\"older's inequality, with $1\leq r\leq p$, and $\frac{1}{r}=\frac{1}{p}+\frac{1}{a}$. Then, using H\"older's interpolation inequality and Lemma \ref{lemma3}, we have
\begin{eqnarray*}
&{}&\left\|\int_0^te^{-(t-s)(-\Delta)^{\beta/2}} f(u(s))\,ds\right\|_{L^p(\mathbb{R}^n)}\\
&{}&\leq C\sum_{k=0}^{\infty}\frac{\lambda^k}{k!}\int_0^t(t-s)^{-\frac{n}{\beta}\left(\frac{1}{r}-\frac{1}{p}\right)}\|u(s)\|_{L^p}\|u(s)\|^{(kp+m-1)\theta}_{L^{q}}\|u(s)\|^{(kp+m-1)(1-\theta)}_{L^{\rho}}\,ds\\
&{}& \leq C\sum_{k=0}^{\infty}\frac{\lambda^k}{k!}\left(\Gamma\left(\frac{\rho}{p}+1\right)\right)^{\frac{(kp+m-1)(1-\theta)}{\rho}}\\
&{}&\qquad \times \int_0^t(t-s)^{-\frac{n}{\beta}\left(\frac{1}{r}-\frac{1}{p}\right)}\|u(s)\|_{\textnormal{exp~$L^p$}}\|u(s)\|^{(kp+m-1)\theta}_{L^{q}}\|u(s)\|^{(kp+m-1)(1-\theta)}_{\textnormal{exp~$L^p$}}\,ds\\
&{}& \leq C\sum_{k=0}^{\infty}\frac{\lambda^k}{k!}\left(\Gamma\left(\frac{\rho}{p}+1\right)\right)^{\frac{(kp+m-1)(1-\theta)}{\rho}}\\
&{}&\times \int_0^t(t-s)^{-\frac{n}{\beta}\left(\frac{1}{r}-\frac{1}{p}\right)}s^{-\sigma(kp+m-1)\theta}\|u(s)\|_{\textnormal{exp~$L^p$}}\left(s^{\sigma}\|u(s)\|\right)^{(kp+m-1)\theta}_{L^{q}}\|u(s)\|^{(kp+m-1)(1-\theta)}_{\textnormal{exp~$L^p$}}\,ds,
\end{eqnarray*}
where
$$\frac{1}{a(kp+m-1)}=\frac{\theta}{q}+\frac{1-\theta}{\rho},\quad 0\leq \theta\leq 1,\quad\hbox{and}\quad p\leq\rho<\infty.$$
By using the fact that $u\in B_\varepsilon$, we get
\begin{eqnarray}\label{eq15}
&{}&\left\|\int_0^te^{-(t-s)(-\Delta)^{\beta/2}} f(u(s))\,ds\right\|_{L^p(\mathbb{R}^n)}\nonumber\\
&{}& \leq C\sum_{k=0}^{\infty}\frac{\lambda^k}{k!}\left(\Gamma\left(\frac{\rho}{p}+1\right)\right)^{\frac{(kp+m-1)(1-\theta)}{\rho}} (M\varepsilon)^{kp+m}\int_0^t(t-s)^{-\frac{n}{\beta}\left(\frac{1}{r}-\frac{1}{p}\right)}s^{-\sigma(kp+m-1)\theta}\,ds\nonumber\\
&{}&=C\sum_{k=0}^{\infty}\frac{\lambda^k}{k!}\left(\Gamma\left(\frac{\rho}{p}+1\right)\right)^{\frac{(kp+m-1)(1-\theta)}{\rho}} (M\varepsilon)^{kp+m}t^{1-\frac{n}{\beta}\left(\frac{1}{r}-\frac{1}{p}\right)-\sigma(kp+m-1)\theta}\nonumber\\
&{}&\qquad \times \int_0^1(1-s)^{-\frac{n}{\beta}\left(\frac{1}{r}-\frac{1}{p}\right)}s^{-\sigma(kp+m-1)\theta}\,ds\nonumber\\
&{}&=C\sum_{k=0}^{\infty}\frac{\lambda^k}{k!}\left(\Gamma\left(\frac{\rho}{p}+1\right)\right)^{\frac{(kp+m-1)(1-\theta)}{\rho}} \nonumber\\
&{}&\qquad\times(M\varepsilon)^{kp+m}\mathcal{B}\left(1-\frac{n}{\beta}\left(\frac{1}{r}-\frac{1}{p}\right);1-\sigma(kp+m-1)\theta\right),\qquad\qquad
\end{eqnarray}
where $\mathcal{B}$ is the beta function, under the following conditions:
$$\frac{n}{\beta}\left(\frac{1}{r}-\frac{1}{p}\right)<1,\qquad \sigma(kp+m-1)\theta<1,\quad\hbox{and}\quad 1-\frac{n}{\beta}\left(\frac{1}{r}-\frac{1}{p}\right)-\sigma(kp+m-1)\theta=0.$$
It remains to prove the existence of $\theta=\theta_k$, $\rho=\rho_k$, $k\geq0$, and $a$. As $q>\frac{(m-1)p}{p-1}$, one can choose 
$$\frac{1-\frac{n(p-1)}{p\beta}}{\sigma(pk+m-1)}<\theta_k<\frac{1}{pk+m-1}\min(m-1,\frac{1}{\sigma}),$$
 and as $\sigma=\frac{1}{m-1}-\frac{n}{\beta q}<\frac{1}{m-1}$; it follows that $\theta_k$ is chosen by
$$\frac{1-\frac{n(p-1)}{p\beta}}{\sigma(pk+m-1)}<\theta_k<\frac{m-1}{pk+m-1}.$$ 
We note that the lower bound of $\theta_k$ is just to be compatible with the condition that $r>1$. For the choice of $\rho_k$, we explain slightly the steps; we need the condition $1-\frac{n}{\beta}\left(\frac{1}{r}-\frac{1}{p}\right)-\sigma(kp+m-1)\theta_k=0$, and as $\frac{1}{r}=\frac{1}{p}+\frac{1}{a}$, so $1-\frac{n}{a\beta}-\sigma(kp+m-1)\theta_k=0$. Then, using the fact that $\frac{1}{a(kp+m-1)}=\frac{\theta_k}{q}+\frac{1-\theta_k}{\rho_k}$ and $\sigma=\frac{1}{m-1}-\frac{n}{\beta q}$, we conclude that $\rho=\rho_k$ is chosen such that
$$\frac{1-\theta_k}{\rho_k}=\frac{\beta}{n(kp+m-1)}-\frac{\beta\theta_k}{n(m-1)}.$$
We note that $\frac{1-\theta_k}{\rho_k}\leq\frac{\beta}{n(m-1)}-\frac{\beta\theta_k}{n(m-1)}=\frac{\beta(1-\theta_k)}{n(m-1)}$ which implies that $\rho_k\geq \frac{n(m-1)}{\beta}\geq p$. Finally, we choose $a>0$ such that
$$\frac{1}{a(kp+m-1)}=\frac{\theta_k}{q}+\frac{1-\theta_k}{\rho_k}.$$
Moreover, for these choice of parameters, 
\begin{equation}\label{eq16}
\mathcal{B}\left(1-\frac{n}{\beta}\left(\frac{1}{r}-\frac{1}{p}\right);1-\sigma(kp+m-1)\theta\right)=\frac{\Gamma\left(1-\frac{n}{\beta}\left(\frac{1}{r}-\frac{1}{p}\right)\right)\Gamma\left(\frac{n}{\beta}\left(\frac{1}{r}-\frac{1}{p}\right)\right)}{\Gamma(1)}\leq C,
\end{equation}
where we have used the fact that $\mathcal{B}(x,y)=\frac{\Gamma(x)\Gamma(y)}{\Gamma(x+y)}$, for every $x,y>0$. We notice also that
$$\theta_k\longrightarrow 0,\,\, \rho_k\longrightarrow \infty\quad\hbox{as}\quad k\rightarrow\infty,$$
then
$$\frac{(kp+m-1)(1-\theta_k)}{p\rho_k}(1+\rho_k)\leq k, \quad \hbox{for all}\,\, k\geq1,$$
this implies, together with the property $\Gamma(x+1)\leq C\,x^{x+\frac{1}{2}}$, $\hbox{for all}\,\, x\geq 1$, that
\begin{equation}\label{eq17}
\left(\Gamma\left(\frac{\rho_k}{p}+1\right)\right)^{\frac{(kp+m-1)(1-\theta_k)}{\rho_k}} \leq C^k\,k!.
\end{equation}
Combining (\ref{eq15}), (\ref{eq16}) and (\ref{eq17}), we obtain
$$
\left\|\int_0^te^{-(t-s)(-\Delta)^{\beta/2}} f(u(s))\,ds\right\|_{L^p(\mathbb{R}^n)} \leq C\sum_{k=0}^{\infty}(C\,\lambda)^k (M\varepsilon)^{kp+m}\leq C(M\varepsilon)^m,
$$
for $\varepsilon$ small enough. This proves (\ref{eq113}). Next, we prove (\ref{eq14}). Using the fact that $\beta=\frac{n(p-1)}{p}$ and Lemma \ref{lemma5}, we have
 $$\left\|\int_0^te^{-(t-s)(-\Delta)^{\beta/2}} f(u(s))\,ds\right\|_{L^\infty(0,\infty;L^\phi(\mathbb{R}^n))}\leq C\|f(u(s))\|_{L^\infty(0,\infty;L^1(\mathbb{R}^n)\cap L^{2p}(\mathbb{R}^n)\cap L^{\frac{2p}{p-1}}(\mathbb{R}^n))}.$$ 
As
 $$|f(u)|\leq C|u|^{m}e^{\lambda |u|^p}= C|u|^{m}\left(e^{\lambda |u|^p}-1\right)+C|u|^{m},$$
so, using $m\geq p$ and a similar calculation as in the case of $\beta<\frac{n(p-1)}{p}$ (see (\ref{eq19})), we conclude that
$$\left\|f(u(t))\right\|_{L^r(\mathbb{R}^n)}\leq C(M\varepsilon)^m,$$
for $r=1,2p,\frac{2p}{p-1}\geq 1$, and all $t>0$. This proves (\ref{eq14}).\\
To estimate the second term in $\Phi(u)$ on $\textnormal{exp~$L^p(\mathbb{R}^n)$}$ in the case of $\beta>\frac{n(p-1)}{p}$, let $b>0$ be the positive number satisfying $b=2\ln(b+1)$, then we can check that
\begin{equation}\label{eq18}
\left(\ln\left((t-s)^{-n/\beta}+1\right)\right)^{-1/p}\leq 2^{1/p}(t-s)^{n/\beta p},\qquad\hbox{for}\,\,0\leq s\leq t-b^{-\beta/n}.
\end{equation}
If $t\leq b^{-\beta/n}$, similarly to (\ref{eq11}), we have
\begin{eqnarray*}
&{}&\left\|\int_0^te^{-(t-s)(-\Delta)^{\beta/2}} f(u(s))\,ds\right\|_{\textnormal{exp~$L^p$}}\\
&{}&\,\leq \int_0^t\left\|e^{-(t-s)(-\Delta)^{\beta/2}} f(u(s))\right\|_{\textnormal{exp~$L^p$}}\,ds\\
&{}&\,\leq  \int_0^t\left(C(t-s)^{-\frac{n}{\beta r}}+1\right)\left(\|f(u(s))\|_r+\|f(u(s))\|_1\right)\,ds,
\end{eqnarray*}
for any $r\geq 1$. Let $r=\frac{p}{p-1}>1$, we get
\begin{eqnarray}\label{eq20}
&{}&\left\|\int_0^te^{-(t-s)(-\Delta)^{\beta/2}} f(u(s))\,ds\right\|_{\textnormal{exp~$L^p$}}\nonumber\\
&{}&\,\leq  \int_0^t\left(C(t-s)^{-\frac{n(p-1)}{\beta p}}+1\right)\left(\|f(u(s))\|_{\frac{p}{p-1}}+\|f(u(s))\|_1\right)\,ds\nonumber\\
&{}&\,\leq \|f(u)\|_{L^\infty(0,\infty;L^1(\mathbb{R}^n)\cap L^{\frac{p}{p-1}}(\mathbb{R}^n))}  \int_0^t\left(C(t-s)^{-\frac{n(p-1)}{\beta p}}+1\right)\,ds\nonumber\\
&{}&\,= \|f(u)\|_{L^\infty(0,\infty;L^1(\mathbb{R}^n)\cap L^{\frac{p}{p-1}}(\mathbb{R}^n))}  \int_0^t\left(Cs^{-\frac{n(p-1)}{\beta p}}+1\right)\,ds\nonumber\\
&{}&\,\leq \|f(u)\|_{L^\infty(0,\infty;L^1(\mathbb{R}^n)\cap L^{\frac{p}{p-1}}(\mathbb{R}^n))}  \int_0^{b^{-\beta/n}}\left(Cs^{-\frac{n(p-1)}{\beta p}}+1\right)\,ds\nonumber\\
&{}&\,=C \|f(u)\|_{L^\infty(0,\infty;L^1(\mathbb{R}^n)\cap L^{\frac{p}{p-1}}(\mathbb{R}^n))},
\end{eqnarray}
where we have used the fact that $\beta>\frac{n(p-1)}{p}$. Then, using $m> p$ and similarly to (\ref{eq19}), we conclude that
$$\left\|f(u(t))\right\|_{L^r(\mathbb{R}^n)}\leq C(M\varepsilon)^m,$$
for $r=1,\frac{p}{p-1}\geq 1$, i.e.
$$\left\|\int_0^te^{-(t-s)(-\Delta)^{\beta/2}} f(u(s))\,ds\right\|_{\textnormal{exp~$L^p$}}=O(\varepsilon).$$
If $t> b^{-\beta/n}$, we have
\begin{eqnarray*}
\left\|\int_0^te^{-(t-s)(-\Delta)^{\beta/2}} f(u(s))\,ds\right\|_{\textnormal{exp~$L^p$}}&\leq& \int_0^t\left\|e^{-(t-s)(-\Delta)^{\beta/2}} f(u(s))\right\|_{\textnormal{exp~$L^p$}}\,ds\\
&=&\int_0^{t-b^{-\beta/n}}\left\|e^{-(t-s)(-\Delta)^{\beta/2}} f(u(s))\right\|_{\textnormal{exp~$L^p$}}\,ds\\
&{}&\,+\int_{t-b^{-\beta/n}}^t\left\|e^{-(t-s)(-\Delta)^{\beta/2}} f(u(s))\right\|_{\textnormal{exp~$L^p$}}\,ds\\
&=:& I+II.
\end{eqnarray*}
Similarly to (\ref{eq20}), using $\beta>\frac{n(p-1)}{p}$ and $m> p$, we have
\begin{eqnarray*}
II&\leq&  \|f(u)\|_{L^\infty(0,\infty;L^1\cap L^{\frac{p}{p-1}})}  \int_{t-b^{-\beta/n}}^t\left(C(t-s)^{-\frac{n(p-1)}{\beta p}}+1\right)\,ds\\
&=&\|f(u)\|_{L^\infty(0,\infty;L^1\cap L^{\frac{p}{p-1}})}  \int_0^{b^{-\beta/n}}\left(Cs^{-\frac{n(p-1)}{\beta p}}+1\right)\,ds \leq C(M\varepsilon)^m.
\end{eqnarray*}
On the other hand, using Proposition \ref{prop1} (ii) and (\ref{eq18}), we have
\begin{eqnarray*}
I&\leq&C\int_0^{t-b^{-\beta/n}}(t-s)^{-\frac{n}{\beta a}}\left(\ln\left((t-s)^{-n/\beta}+1\right)\right)^{-1/p}\left\|f(u(s))\right\|_{L^a}\,ds\\
&\leq&C\int_0^{t}(t-s)^{-\frac{n}{\beta}(\frac{1}{a}-\frac{1}{p})}\left\|f(u(s))\right\|_{L^a}\,ds,
\end{eqnarray*}
where $1\leq a\leq p$. Apply the same calculation done above to obtain (\ref{eq113}) (with same conditions), we conclude that
$$I=O(\varepsilon).$$
This implies that
$$\left\|\int_0^te^{-(t-s)(-\Delta)^{\beta/2}} f(u(s))\,ds\right\|_{\textnormal{exp~$L^p$}}=O(\varepsilon),$$
in the case of $t> b^{-\beta/n}$, therefore
$$\left\|\int_0^te^{-(t-s)(-\Delta)^{\beta/2}} f(u(s))\,ds\right\|_{L^\infty(0,\infty;\textnormal{exp~$L^p$})}=O(\varepsilon),\quad\hbox{for all}\,\,t>0.$$
It remains to prove that
$$t^\sigma\left\|\int_0^te^{-(t-s)(-\Delta)^{\beta/2}} f(u(s))\,ds\right\|_{L^q(\mathbb{R}^n)}=O(\varepsilon),$$
for every $t>0$, to conclude that $\Phi(u)\in B_\varepsilon$. This follows similarly as in (\ref{eq21}) below by using the fact that $f(0)=0$.\\

\noindent$\bullet$ {\bf $\Phi$ is a contraction.} Let $u,v\in B_\varepsilon$. By (\ref{P_2}), we obtain
$$
t^\sigma\left\|\Phi(u)-\Phi(v)\right\|_{L^q(\mathbb{R}^n)}\leq C\,t^\sigma\int_0^t(t-s)^{-\frac{n}{\beta}\left(\frac{1}{r}-\frac{1}{q}\right)}\left\|f(u(s))-f(v(s))\right\|_r\,ds,
$$
for every $1\leq r\leq q$. From our assumption (\ref{eq4}), we have
\begin{eqnarray*}
|f(u)-f(v)|&\leq & C|u-v|(|u|^{m-1}e^{\lambda |u|^p}+|v|^{m-1}e^{\lambda |v|^p})\\
&=& C|u-v|\left(|u|^{m-1}\sum_{k=0}^{\infty}\frac{\lambda^k}{k!}|u|^{kp}+|v|^{m-1}\sum_{k=0}^{\infty}\frac{\lambda^k}{k!}|v|^{kp}\right)\\
 &=& C\sum_{k=0}^{\infty}\frac{\lambda^k}{k!}|u-v|\left(|u|^{kp+m-1}+|v|^{kp+m-1}\right).
 \end{eqnarray*}
Using H\"older's inequality and H\"older's interpolation inequality, we get 
\begin{eqnarray*}
\|f(u)-f(v)\|_r&\leq&  C\sum_{k=0}^{\infty}\frac{\lambda^k}{k!}\left\|(u-v)\left(|u|^{kp+m-1}+|v|^{kp+m-1}\right)\right\|_r\\
&\leq&  C\sum_{k=0}^{\infty}\frac{\lambda^k}{k!}\|u-v\|_q\left(\|u\|_{a(kp+m-1)}^{kp+m-1}+\|v\|_{a(kp+m-1)}^{kp+m-1}\right)\\
&\leq&  C\sum_{k=0}^{\infty}\frac{\lambda^k}{k!}\|u-v\|_q\left(\|u\|_{q}^{(kp+m-1)\theta}\|u\|_{\rho}^{(kp+m-1)(1-\theta)}\right.\\
&{}&\left.\qquad\qquad\qquad\qquad+\|v\|_{q}^{(kp+m-1)\theta}\|v\|_{\rho}^{(kp+m-1)(1-\theta)}\right),
 \end{eqnarray*}
where
$$\frac{1}{r}=\frac{1}{q}+\frac{1}{a}\quad \hbox{and}\quad \frac{1}{a(kp+m-1)}=\frac{\theta}{q}+\frac{1-\theta}{\rho},\quad \hbox{for all}\,\, 0\leq \theta\leq 1.$$
Using Lemma \ref{lemma3}, assuming that $p\leq\rho<\infty$, we infer that
\begin{eqnarray*}
&{}&\|f(u)-f(v)\|_r\\
&{}&\leq  C\sum_{k=0}^{\infty}\frac{\lambda^k}{k!}\left(\Gamma\left(\frac{\rho}{p}+1\right)\right)^{\frac{(kp+m-1)(1-\theta)}{\rho}}\\
&{}&\quad\times\|u-v\|_q\left(\|u\|_{q}^{(kp+m-1)\theta}\|u\|_{\textnormal{exp~$L^p$}}^{(kp+m-1)(1-\theta)}+\|v\|_{q}^{(kp+m-1)\theta}\|v\|_{\textnormal{exp~$L^p$}}^{(kp+m-1)(1-\theta)}\right).
 \end{eqnarray*}
So
\begin{eqnarray*}
&{}&t^\sigma\left\|\Phi(u)-\Phi(v)\right\|_{L^q(\mathbb{R}^n)}\\
&{}&\leq  C\sum_{k=0}^{\infty}\frac{\lambda^k}{k!}\left(\Gamma\left(\frac{\rho}{p}+1\right)\right)^{\frac{(kp+m-1)(1-\theta)}{\rho}}\\
&{}&\,\,\times t^\sigma\int_0^t(t-s)^{-\frac{n}{\beta}\left(\frac{1}{r}-\frac{1}{q}\right)}s^{-\sigma}s^\sigma\|u-v\|_qs^{-\sigma(kp+m-1)\theta}\\
&{}&\,\,\times\left(\left(s^\sigma\|u\|_{q}\right)^{(kp+m-1)\theta}\|u\|_{\textnormal{exp~$L^p$}}^{(kp+m-1)(1-\theta)}+\left(s^\sigma\|v\|_{q}\right)^{(kp+m-1)\theta}\|v\|_{\textnormal{exp~$L^p$}}^{(kp+m-1)(1-\theta)}\right)\,ds\\
&{}&\leq Cd(u,v)(\varepsilon M)^{m-1}\sum_{k=0}^{\infty}(\varepsilon M)^{kp}\frac{\lambda^k}{k!}\left(\Gamma\left(\frac{\rho}{p}+1\right)\right)^{\frac{(kp+m-1)(1-\theta)}{\rho}}\\
&{}&\,\,\times\mathcal{B}\left(1-\frac{n}{\beta}\left(\frac{1}{r}-\frac{1}{q}\right);1-\sigma(1+(kp+m-1)\theta)\right),
 \end{eqnarray*}
 where we have used the fact that $u,v\in B_\varepsilon$, under the following conditions:
 $$1-\frac{n}{\beta}\left(\frac{1}{r}-\frac{1}{q}\right)-\sigma(kp+m-1)\theta=0,\quad \frac{n}{\beta}\left(\frac{1}{r}-\frac{1}{q}\right)<1,\quad \hbox{and}\quad \sigma(1+(kp+m-1)\theta)<1.$$
As above, for all $k\geq 0$, we choose first $\theta=\theta_k\geq0$ such that
$$\frac{1-\frac{n(q-1)}{q\beta}}{\sigma(pk+m-1)}<\theta<\frac{1}{pk+m-1}\min(m-1,\frac{1-\sigma}{\sigma}),$$
where we have used the fact that $q>\frac{(m-1)p}{p-1}\geq m$. Next, we choose $\rho=\rho_k$ such that
$$\frac{1-\theta_k}{\rho_k}=\frac{\beta}{n(kp+m-1)}-\frac{\beta\theta_k}{n(m-1)},$$
and finally, we choose $a>0$ such that
$$\frac{1}{a(kp+m-1)}=\frac{\theta_k}{q}+\frac{1-\theta_k}{\rho_k}.$$
To ensure that $\sigma<1$, we also suppose the following condition 
$$q<\frac{n(m-1)}{\beta(2-m)_+},$$
where $(\cdotp)_+$ stands for the positive part. Moreover, for these choice of parameters, 
$$
\mathcal{B}\left(1-\frac{n}{\beta}\left(\frac{1}{r}-\frac{1}{p}\right);1-\sigma(1+(kp+m-1)\theta)\right)=\frac{\Gamma\left(1-\frac{n}{\beta}\left(\frac{1}{r}-\frac{1}{p}\right)\right)\Gamma\left(\frac{n}{\beta}\left(\frac{1}{r}-\frac{1}{p}\right)\right)}{\Gamma(\frac{m-2}{m-1}+\frac{n}{\beta q})}\leq C,
$$
and
$$
\left(\Gamma\left(\frac{\rho_k}{p}+1\right)\right)^{\frac{(kp+m-1)(1-\theta_k)}{\rho_k}} \leq C^k\,k!.
$$
This implies that
\begin{equation}\label{eq21}
t^\sigma\left\|\Phi(u)-\Phi(v)\right\|_{L^q(\mathbb{R}^n)}\leq Cd(u,v)(M\varepsilon)^{m-1}\sum_{k=0}^{\infty}(C\,\lambda)^k (M\varepsilon)^{kp}\leq \frac{1}{2} d(u,v),
\end{equation}
for $\varepsilon$ small enough. This completes the proof the existence of global solution in Theorem \ref{theo2} in the case of $\beta\geq\frac{n(p-1)}{p}$. The estimation (\ref{eq35}) follows from $u\in B_\varepsilon$.  \hfill$\square$\\


\subsection{\bf {Proof of the property (\ref{eq5}) in Theorem \ref{theo2}}}\label{subsec4.3}
We now prove the continuity of solution at zero. Let $q$ be a positive number such that $q> \max\{\frac{n}{\beta},1\}$. From the embedding $L^p(\mathbb{R}^n)\cap L^\infty(\mathbb{R}^n)\hookrightarrow \textnormal{exp}\, L^p(\mathbb{R}^n)$ (Lemma \ref{lemma1}), and $L^p-L^p$, $L^q-L^\infty$ estimates (\ref{P_2}), we have
\begin{eqnarray}\label{eq22}
&{}&\|u(t)-e^{-t(-\Delta)^{\beta/2}}u_0\|_{\textnormal{exp~$L^p$}}\nonumber\\
&{}&\leq\displaystyle  \int_0^t \left\|e^{-(t-s)(-\Delta)^{\beta/2}}f(u(s))\right\|_{\textnormal{exp~$L^p$}}\,ds\nonumber\\
&{}&\leq C\displaystyle  \int_0^t \left\|e^{-(t-s)(-\Delta)^{\beta/2}}f(u(s))\right\|_{L^p}\,ds+C\displaystyle  \int_0^t \left\|e^{-(t-s)(-\Delta)^{\beta/2}}f(u(s))\right\|_{L^\infty}\,ds\nonumber\\
&{}&\leq C\displaystyle  \int_0^t \left\|f(u(s))\right\|_{L^p}\,ds+C\displaystyle  \int_0^t (t-s)^{-\frac{n}{\beta q}}\left\|f(u(s))\right\|_{L^q}\,ds.
\end{eqnarray}
Let us estimate $\|f(u)\|_{L^r}$, for $r=p,q\geq 1$. We have
$$|f(u)|\leq C|u|^{m}e^{\lambda |u|^p}= C|u|^{m}\left(e^{\lambda |u|^p}-1\right)+C|u|^{m},$$
then, by H\"older's inequality, we obtain
\begin{eqnarray*}
\left\|f(u)\right\|_{L^r(\mathbb{R}^n)}&\leq& C \left\|u\right\|^{m}_{L^{2mr}(\mathbb{R}^n)} \left\|e^{\lambda |u|^p}-1\right\|_{L^{2r}(\mathbb{R}^n)}+C\left\|u\right\|^{m}_{L^{mr}(\mathbb{R}^n)}\\
&\leq& C \left\|u\right\|^{m}_{\textnormal{exp~$L^p(\mathbb{R}^n)$}} \left\|e^{\lambda |u|^p}-1\right\|_{L^{2r}(\mathbb{R}^n)}+C\left\|u\right\|^{m}_{\textnormal{exp~$L^p(\mathbb{R}^n)$}},
\end{eqnarray*}
where we have used Lemma \ref{lemma3} and $2mr\geq mr\geq m\geq p$. Next, using Lemma \ref{lemma6} and the fact that $u\in E_\varepsilon$ (or $u\in B_\varepsilon$), we have
\begin{equation}\label{eq23}
\left\|f(u)\right\|_{L^r(\mathbb{R}^n)}\leq C \left\|u\right\|^{m}_{\textnormal{exp~$L^p(\mathbb{R}^n)$}}(1+2C\lambda r(\varepsilon)^p)^{1/2r}\leq C\left\|u\right\|^{m}_{\textnormal{exp~$L^p(\mathbb{R}^n)$}}.
\end{equation}
Substituting (\ref{eq23}) in (\ref{eq22}), we obtain
 \begin{eqnarray*}
\|u(t)-e^{-t(-\Delta)^{\beta/2}}u_0\|_{\textnormal{exp~$L^p$}}&\leq&C\displaystyle  \int_0^t \left\|u\right\|^{m}_{\textnormal{exp~$L^p$}}\,ds+C\displaystyle  \int_0^t (t-s)^{-\frac{n}{\beta q}}\left\|u\right\|^{m}_{\textnormal{exp~$L^p$}}\,ds\\
&\leq&C\,t \left\|u\right\|^{m}_{L^\infty(0,\infty;\textnormal{exp~$L^p$})}+C\,t^{1-\frac{n}{\beta q}}\left\|u\right\|^{m}_{L^\infty(0,\infty;\textnormal{exp~$L^p$})}\\
&\leq&C\,t +C\,t^{1-\frac{n}{\beta q}}\longrightarrow 0\qquad\hbox{as}\quad t\rightarrow0.
\end{eqnarray*}
This completes the proof of (\ref{eq5}). \hfill$\square$\\


\subsection{\bf {Proof of the weak$^{*}$ convergence in Theorem \ref{theo2}}}\label{subsec4.4}
We complete the proof of Theorem \ref{theo2} by showing the continuity at $t=0$ in the weak$^{*}$ sense. Let $X:=L^1(\ln L)^{1/p}(\mathbb{R}^n)$ be the pre-dual space of $\textnormal{exp~$L^p$}$. It is known that $X$ is a Banach space and $C^\infty_0(\mathbb{R}^n)$ is dense in $X$ (cf. \cite{Adams}). Let $\varphi\in X$. By H\"older's inequality for the Orlicz space, we have
 \begin{eqnarray*}
\left|\int_{\mathbb{R}^n}\left(e^{-t(-\Delta)^{\beta/2}}u_0(x)-u_0(x)\right)\varphi(x)\,dx\right|&=&\left|\int_{\mathbb{R}^n}u_0(x)\left(e^{-t(-\Delta)^{\beta/2}}\varphi(x)-\varphi(x)\right)\,dx\right|\\
&\leq& 2\|u_0\|_{\textnormal{exp~$L^p$}}\left\|e^{-t(-\Delta)^{\beta/2}}\varphi-\varphi\right\|_X.
 \end{eqnarray*}
Since $C^\infty_0(\mathbb{R}^n)$ is dense in $X$, so by applying similar calculations as in the proof of Proposition \ref{prop2}, we conclude that 
$$\lim_{t\rightarrow0}\left\|e^{-t(-\Delta)^{\beta/2}}\varphi-\varphi\right\|_X=0.$$
This completes the weak$^{*}$ convergence.\hfill$\square$\\

 \section*{Acknowledgements}
The authors wish to thank the anonymous referee for his/her valuable comments which helped to improve the article.

\medskip
Received xxxx 20xx; revised xxxx 20xx.
\medskip

\end{document}